\documentclass{amsart}
\usepackage{amsmath,amsthm}
\usepackage{amsfonts,amssymb}
\usepackage{accents}
\usepackage{enumerate}
\usepackage{accents,color}
\usepackage{graphicx}

\newcommand{\lvt}{\left|\kern-1.35pt\left|\kern-1.3pt\left|}
\newcommand{\rvt}{\right|\kern-1.3pt\right|\kern-1.35pt\right|}

\newtheorem{thm}{Theorem}[section]
\newtheorem{cor}[thm]{Corollary}
\newtheorem{lem}[thm]{Lemma}
\newtheorem{prop}[thm]{Proposition}

\newtheorem{defn}[thm]{Definition}

\theoremstyle{remark}
\newtheorem{rem}{Remark}[section]

 \def\la{{\langle}}
 \def\ra{{\rangle}}
 \def\ve{{\varepsilon}}

 \def\d{\mathrm{d}}
 
 \def\sph{{\mathbb{S}^{d-1}}}

 \def\sL{{\mathsf L}}
 
 \def\sP{{\mathsf P}}
 
 \def\sS{{\mathsf S}}
 
 \def\bs{{\mathsf b}}
 \def\sc{{\mathsf c}}
 \def\sd{{\mathsf d}}
 \def\sm{{\mathsf m}}
 \def\sw{{\mathsf w}}

 \def\fD{{\mathfrak D}}
 
 \def\a{{\alpha}}
 \def\b{{\beta}}
 \def\g{{\gamma}}
 \def\k{{\kappa}}
 \def\t{{\theta}}
 \def\l{{\lambda}}
 \def\o{{\omega}}
 \def\s{\sigma}
 \def\la{{\langle}}
 \def\ra{{\rangle}}
 \def\ve{{\varepsilon}}

 \def\kb{{\mathbf k}}

 \def\vb{{\mathbf v}}

 \def\Jb{{\mathbf J}}
 
 \def\Lb{{\mathbf L}}
 \def\Pb{{\mathbf P}}

 \def\CD{{\mathcal D}}
 
 \def\CH{{\mathcal H}}

 \def\CV{{\mathcal V}}

 \def\BB{{\mathbb B}}

 \def\NN{{\mathbb N}}

 \def\RR{{\mathbb R}}
 \def\SS{{\mathbb S}}
 \def\TT{{\mathbb T}}
 \def\VV{{\mathbb V}}

      \def\proj{\operatorname{proj}}

\def\lla{\langle{\kern-2.5pt}\langle}      
\def\rra{\rangle{\kern-2.5pt}\rangle}

\newcommand{\wh}{\widehat}

\def\fD{{\mathfrak D}}

\def\f{\frac}

\graphicspath{{./}}
\begin{document}

\title[Bernstein inequality on conic domains and triangles]
{Bernstein inequality on conic domains and triangles}

\author{Yuan~Xu}
\address{Department of Mathematics, University of Oregon, Eugene, 
OR 97403--1222, USA}
\email{yuan@uoregon.edu} 
\thanks{The author is partially supported by Simons Foundation Grant \#849676 and by an Alexander von Humboldt 
award.}

\date{\today}  

\begin{abstract} 
We establish weighted Bernstein inequalities in $L^p$ space for the doubling weight on the conic surface 
$\VV_0^{d+1} = \{(x,t): \|x\| = t, x \in \RR^d, t\in [0,1]\}$ as well as on the solid cone bounded by the conic 
surface and the hyperplane $t =1$, which becomes a triangle on the plane when $d=1$. While the inequalities 
for the derivatives in the $t$ variable behave as expected, there are inequalities for the derivatives in the 
$x$ variables that are stronger than what one may have expected. As an example, on the triangle 
$\{(x_1,x_2): x_1 \ge 0, \, x_2 \ge 0,\, x_1+x_2 \le 1\}$,  the usual Bernstein inequality for the derivative
$\partial_1$ states that $\|\phi_1 \partial_1 f\|_{p,w} \le  c n \|f\|_{p,w}$ with $\phi_1(x_1,x_2):= x_1(1-x_1-x_2)$, 
whereas our new result gives 
$$\| (1-x_2)^{-1/2} \phi_1 \partial_1 f\|_{p,w} \le  c n \|f\|_{p,w}.$$
 The new inequality is stronger and points out a phenomenon unobserved hitherto for polygonal domains. 
\end{abstract}
\maketitle

\section{Introduction} 
\setcounter{equation}{0}

The Bernstein inequalities are fundamental in approximation theory, as seen in the inverse estimate in the 
characterization of best approximation and numerous other applications (see, for example, \cite{DG, DT}). 
Starting from algebraic polynomials on $[0,1]$, the Bernstein or Bernstein-Markov inequalities have
been refined and generalized extensively by many authors. Most notable recent extensions are the inequalities
in $L^p$ norm with the doubling weight, initiated in \cite{MT1}, and inequalities for multivariable polynomials
on various domains, such as polytopes, convex domains, and domains with smooth boundaries; see 
\cite{B1, B2, BLMT, BM, BLMR, Dai1, DP2, DaiX, Kroo0, Kroo1, Kroo2, Kroo3, KR, LWW, T1, T2, X05} 
and reference therein. The purpose of this paper is to establish the Bernstein inequalities, in uniform norm 
and $L^p$ norm with doubling weight, for polynomials on the conic surface
$$
  \VV_0^{d+1} = \left \{(x,t): \|x\| = t, \quad x \in \RR^d, \, 0 \le t \le 1 \right\},
$$ 
and on the solid cone $\VV^{d+1}$ bounded by $\VV_0^{d+1}$ and the hyperplane $t =1$, as well as on 
the planar triangle. The domain has a singularity at the vertex. How singularity affects the Bernstein inequalities 
motivates our study. 

We will establish several inequalities that demonstrate the impact of the singularity at the vertex for the conic 
domains. Our results also lead to new inequalities on triangle domains that are stronger than those known in 
the literature, which are somewhat unexpected and reveal new phenomena hitherto unnoticed. To be more 
precise, let us first recall the Bernstein inequalities in $L^p$ norm with doubling weight on the interval \cite{MT1}
and on the unit sphere \cite{Dai1, DaiX}. Let $\ell$ be a positive integer. For a doubling weight $w$ on $[0,1]$ 
and $1 \le p < \infty$, it is known \cite[Theorem 7.3]{MT1} that 
\begin{equation}\label{eq:Bernstein[0,1]}
  \left \|f^{(\ell)} \right\|_{L^p([0,1],w)} \le c n^{2 \ell} \|f\|_{L^p([0,1],w)}, \qquad  \deg f \le n;
\end{equation}
moreover, setting $\varphi(x) = \sqrt{x(1-x)}$, then \cite[Theorem 7.4]{MT1} 
\begin{equation}\label{eq:Bernstein[0,1]phi}
   \left \|\varphi^\ell f^{(\ell)}  \right\|_{L^p([0,1],w)} \le c n^\ell \|f\|_{L^p([0,1],w)}, \qquad  \deg f \le n.
\end{equation}
For polynomials on the unit sphere $\sph$ and a doubling weight $w$ on $\sph$, the Bernstein inequalities are
of the form \cite[Theorem 5.5.2]{DaiX} 
\begin{equation}\label{eq:BernsteinSph}
  \left  \|D_{i,j}^\ell f  \right \|_{L^p(\sph,w)} \le c n^\ell \|f\|_{L^p(\sph,w)},  \qquad  \deg f \le n, 
\end{equation}
where the derivatives $D_{i,j}$ are the angular derivatives on the unit sphere defined by 
\begin{equation}\label{eq:Dij}
   D_{i,j} = x_i \f{\partial}{\partial x_j} - x_j \f{\partial}{\partial x_i}, \quad 1 \le i < j \le d.
\end{equation}
These inequalities also hold for the uniform norm in the respective domain. 

The Bernstein inequalities on the conic surface are closely related to the above inequalities. Indeed, 
parametrizing $\VV_0^{d+1}$ as $\{(t \xi, t): \xi \in \sph, 0 \le t \le 1\}$, we see that the space of polynomials 
on $\VV_0^{d+1}$ contains the subspace of polynomials in the $t$ variable as well as spherical polynomials 
on the unit sphere $\sph$. Moreover, writing $(x,t) = (t\xi, t) \in \VV_0^{d+1}$ with $\xi \in \sph$, it follows that 
the derivatives on the surface $\VV_0^{d+1}$ are the partial derivative in the $t$ variable and the angular 
derivatives $D_{i,j}$ in the $x$ variable. The conic surface, however, is not a smooth surface because of its 
vertex at the origin of $\RR^{d+1}$. Our result shows that the inequality for the partial derivative in the $t$ 
variable behaves as expected, but the Bernstein inequality for $D_{i,j}$ on $\VV_0^{d+1}$ turns out to satisfy 
$$
  \left \|\frac{1}{\sqrt{t}^\ell} D_{i,j}^\ell f \right \|_{L^p(\VV_0^{d+1},\sw)} \le c n^\ell \|f\|_{L^p(\VV_0^{d+1}, \sw)}, 
       \quad \ell = 1,2, \quad  \deg f \le n,
$$
where $\sw$ is a doubling weight on the conic surface, but not for $\ell > 2$ in general. The factor $1/\sqrt{t}$ 
reflects the geometry of the conic surface. A similar result will also be established on the solid cone $\VV^{d+1}$.
For $d=1$, the cone $\VV^2$ is a triangle and it is mapped, by an affine transform, to the standard triangle 
$$
\TT^2 = \{(y_1,y_2): y_1\ge 0, y_2\ge 0, y_1+y_2 \le 1\}
$$ 
in $\RR^2$, so that we obtain new Bernstein inequalities on the triangle. In the literature, a typical Bernstein inequality 
on the triangle takes the form
\begin{equation}\label{eq:triangleS1}
\|\phi_{1} \partial_1 f \|_{L^p(\TT^2, w)}\le c \|f\|_{L^p(\TT^2, w)}, \quad \phi_1(y) =  \sqrt{y_1(1-y_1-y_2)},\quad
\deg f \le n,
\end{equation}
for $\partial_1 = \frac{\partial} {\partial y_1}$ for example (cf. \cite{BX, DZ, DT}), and it holds also for the uniform 
norm. The inequality \eqref{eq:triangleS1} is accepted as an appropriate generalization of the Bernstein inequality 
\eqref{eq:Bernstein[0,1]phi} on the triangle, since $\phi_1$ takes into account the boundary of $\TT^2$, just like
$\varphi$ on $[0,1]$ in \eqref{eq:Bernstein[0,1]}. Nevertheless, we obtain a stronger inequality:
\begin{equation}\label{eq:triangleS}
\left \| \frac{1} {\sqrt{1-y_2}}  \phi_{1} \partial_1 f \right\|_{L^p(\TT^2, w)}\le c \|f\|_{L^p(\TT^2, w)}, \quad 
   \deg f \le n.
\end{equation}
The inequality \eqref{eq:triangleS} is somewhat surprising. Indeed, the additional factor $\f{1}{\sqrt{1-y_2}}$ appears 
to be a new phenomenon that has not been observed before. In hindsight, the inequality \eqref{eq:triangleS1}
takes into account the boundary of the triangle, whereas the new ones take into account the singularity at the 
corners as well. 

On several regular domains, including the conic domains, there exist second-order differential operators that have
orthogonal polynomials as eigenfunctions. For the unit sphere, it is the Laplace-Beltrami operator. On the conic 
surface $\VV_0^{d+1}$, it is given by 
$$
   \Delta_{0,\g} =  t(1-t)\partial_t^2 + \big( d-1 - (d+\g)t \big) \partial_t+ t^{-1} \sum_{1\le i < j \le d } \left( D_{i,j}^{(x)} \right)^2,
$$
and it has the orthogonal polynomials with respect to $t^{-1}(1-t)^\g$ on $\VV_0^{d+1}$ as the eigenfunctions. 
The Bernstein inequality for such an operator can be established easily. For example, it was shown in 
\cite[Theorem 3.1.7]{X21} that
$$
 \left \|(- \Delta_{0,\g})^\ell \right \|_{L^p(\VV_0^{d+1},\sw)} \le c n^{2\ell}  \|f\|_{L^p(\VV_0^{d+1}, \sw)}, \quad \ell =1,2,\ldots. 
$$
The proof relies only on the eigenvalues of $\Delta_{0,\g}$ and follows from a general result for all such operators 
on localizable homogeneous spaces. As it is, the dependence on the domain is opaque and hidden in the proof. 
This inequality nevertheless motivates our study on the Bernstein inequalities for the first-order derivatives. It 
can also be used to show the sharpness of our new result for special values of $\ell$ when $p=2$. 

Our analysis relies on orthogonal structures on the conic domains \cite{X20}. It is part of an ongoing program 
that aims at extending the results in approximation theory and harmonics analysis from the unit sphere to 
quadratic surfaces of revolution \cite{X20, X21a, X21, X21b}. Our main tools are the closed-form formula for the 
reproducing kernels in \cite{X20} and the highly localized kernel studied in \cite{X21}. The approach 
follows the analysis on the unit sphere \cite{DaiX}, but the conic domain has its intrinsic complexity. For 
example, the distance function on $\VV_0^{d+1}$ is not intuitively evident and satisfies a formula much 
more involved than the geodesic distance on the unit sphere (see \eqref{eq:distV0} below). The study in 
\cite{X21} establishes a framework, assuming the existence of highly localized kernels, for approximation 
and localized frames on a localizable homogenous space and, along the way, provides a toolbox for carrying 
out analysis on localizable homogeneous spaces. For the conic domains, the highly localized kernels 
are established by delicate estimates. For our proof of the Bernstein inequalities, we shall establish sharp 
estimates for the derivatives of these kernels. 

The paper is organized as follows. The main results are stated and discussed in Section 2. The proofs 
of the main results are in Section 3 for the conic surface and Section 4 for the solid cone.  

Throughout the paper, we let $c, c_1, c_2, \ldots$ as positive constants that depend only on fixed parameters 
and their values could change from line to line. We write $A \sim B$ if $c_1 B \le A \le c_2 B$.

\section{Main results}\label{sec:Main}
\setcounter{equation}{0}

We state and discuss our main results in this section. The Bernstein inequalities on the conic surface 
are stated and discussed in the first subsection and those on the solid cone are stated in the second 
subsection. The inequalities on the triangle follow from those on the solid cone and will be discussed 
in the third subsection.

\subsection{Main results on the conic surface} 
Parametrizing $\VV_0^{d+1}$ as $\{(t \xi, t): \xi \in \sph, 0 \le t \le 1\}$, we see that the first order derivative
on the conic surfaces are $\partial_t = \f{\partial}{\partial t}$ in the $t$ variable and the angular derivatives
$D_{i,j}$, $1\le i, j \le d$, in the $x$ variables, which we denote by $D_{i,j}^{(x)}$. The latter are called angular 
derivatives since if $\t_{i,j}$ is the angle of polar coordinates in the $(x_i,x_j)$-plane, defined by 
$(x_i,x_j) = r_{i,j} (\cos \t_{i,j},\sin \t_{i,j})$, where $r_{i,j} \ge 0$ and $0 \le \t_{i,j} \le 2 \pi$, then a quick 
computation shows that
$$
 D_{i,j} = x_i  \frac{\partial}{\partial x_j} - x_j \frac{\partial}{\partial x_i} = \frac{\partial}{\partial \t_{i,j}},  \qquad 1 \le i,j \le d.
$$
In particular, they are independent of the scale of $x$. We state our Bernstein inequalities for the derivatives 
$\partial_t$ and $D_{i,j}^{(x)}$ on the conic surface. 

Let $\sw$ be a doubling weight on $\VV_0^{d+1}$; see Subsection \ref{sec:OPconicSurface} for the definition. 
We denote by $\|\cdot\|_{p,\sw}$ the weighted $L^p$ norm on 
$\VV_0^{d+1}$
$$
  \|f\|_{p,\sw} = \left( \int_{\VV_0^{d+1}} | f(x,t)|^p \sw(x,t) \d \sm (x,t)\right)^{1/p}, \qquad 1 \le p < \infty, 
$$
where $\sm$ denotes the Lebesgue measure on the conic surface. For $ p =\infty$, we denote by 
$\|\cdot\|_\infty$ the uniform norm on $\VV_0^{d+1}$. For $n = 0,1,2, \ldots$, let $\Pi_n(\VV_0^{d+1})$ be 
the space of polynomials of degree at most $n$ on $\VV_0^{d+1}$.

\begin{thm} \label{thm:BI-V0}
Let $\sw$ be a doubling weight on $\VV_0^{d+1}$ and let $f\in \Pi_n (\VV_0^{d+1})$. For $\ell \in \NN$ and 
$1 \le p < \infty$, 
\begin{equation} \label{eq:BI-V0-1}
     \|\partial_t^\ell  f\|_{p,\sw} \le c_p  n^{2 \ell} \|f\|_{p,\sw},
\end{equation}
and, with $\varphi(t) = \sqrt{t(1-t)}$,  
\begin{equation}\label{eq:BI-V0-2}
     \|\phi^\ell \partial_t^\ell  f\|_{p,\sw} \le c_p  n^{\ell} \|f\|_{p,\sw}.
\end{equation}
Moreover, for $D_{i,j} = D_{i,j}^{(x)}$, $1 \le i < j \le d$, and $\ell \in \NN$, 
\begin{equation}\label{eq:BI-V0-4}
   \left \|D_{i,j}^\ell f\right \|_{p,\sw} \le c_p  n^\ell \|f\|_{p,\sw},  
\end{equation}
and, furthermore, for $\ell =1, 2$,
\begin{equation}\label{eq:BI-V0-3}
   \left \|  \frac{1}{\sqrt{t}^\ell } D_{i,j}^\ell f\right \|_{p,\sw} \le c_p  n^\ell \|f\|_{p,\sw}, 
\end{equation}
but not for $\ell  \ge  3$ in general. Finally, these inequalities also hold when the norm is the 
uniform norm $\|\cdot\|_\infty$ on $\VV_0^{d+1}$. 
\end{thm}
 
The inequality \eqref{eq:BI-V0-3} is stronger than \eqref{eq:BI-V0-4} when $\ell = 1, 2$, since $0 \le t \le 1$. 
Moreover, for $\ell > 2$, the inequality \eqref{eq:BI-V0-4} follows from iteration. It is worth mentioning that 
the inequality \eqref{eq:BI-V0-3} for $\ell =1$ cannot be iterated to obtain the inequality for $\ell =2$, since 
the presence of $1/\sqrt{t}$ means that the function involved in the iteration is no longer a polynomial. 

We now discuss these inequalities in light of the spectral operator $\Delta_{0,\g}$ on the conic surface, which
is defined for $\g > -1$ by 
\begin{equation}\label{eq:LB-operator}
\Delta_{0,\g}:= t(1-t)\partial_t^2 + \big( d-1 - (d+\g)t \big) \partial_t+ t^{-1} \Delta_0^{(\xi)}, 
\end{equation}
where $\xi = x/t$ and $\Delta_0^{(\xi)}$ denote the Laplace-Beltrami operator $\Delta_0$ on the unit sphere
$\sph$ in $\xi$ variables. For $\g > -1$, the operator $- \Delta_{0,\g}$ is non-negative and has orthogonal 
polynomials with respect to the weight function $t^{-1}(1-t)^\g$ on $\VV_0^{d+1}$ as eigenfunctions; see 
Theorem \ref{thm:Jacobi-DE-V0} below. The latter property allows us to define $(-\Delta_{0,r})^{\a}$ for
$\a \in \RR$. It is shown in \cite[Theorem 3.1.7]{X21} that the following Bernstein inequality holds. 

\begin{thm} \label{thm:BernsteinLB}
Let $\sw$ be a doubling weight on $\VV_0^{d+1}$. Let $\g >  -1$. Then, for $r > 0$ and $1 \le p \le \infty$, 
\begin{equation}\label{eq:BernsteinLB}
  \| (-  \Delta_{0,\g})^{\f r 2}f \|_{p,\sw} \le c n^r \|f\|_{p,\sw}, \qquad \forall f\in \Pi_n(\VV_0^{d+1}).
\end{equation}
\end{thm}

The proof of this theorem, however, does not use the explicit formula of the operator $\Delta_{0,\g}$ 
and relies only on the fact that the operator has orthogonal polynomials as eigenfunctions. 
The operator $\Delta_{0,\g}$ in \eqref{eq:LB-operator} contains a factor $t^{-1} \Delta_0^{(\xi)}$. 
It is known that the Laplace-Beltrami operator $\Delta_0$ can be decomposed in terms of angular 
derivatives (cf. \cite[Theorem 1.8.2]{DaiX}), 
\begin{equation}\label{eq:D=Dij} 
   \Delta_0 = \sum_{1\le i < j \le d} D_{i,j}^2.
\end{equation}
Hence, the factor $t^{-1} \Delta_0^{(\xi)}$ in the operator $\Delta_{0,\g}$ can be decomposed as a sum of 
$\frac{1}{\sqrt{t}^{2}} (D_{i,j}^{(\xi)})^2$. Thus, the inequality \eqref{eq:BernsteinLB} with $r=2$ follows
from inequalities in Theorem \ref{thm:BI-V0}. Notice, however, that \eqref{eq:BI-V0-3} holds only for 
$\ell = 1$ and $2$. This makes the higher order Bernstein inequality for the spectral operator 
$\Delta_{0,\g}$ that much more special. 

The operator $- \Delta_{0,\g}$ on the conic surface $\VV_0^{d+1}$ is self-adjoint, which is evident from
the following identity that relies on the first-order derivatives on the conic surface. 

\begin{prop} \label{prop:self-adJacobi}
For $\g > -1$, let $\sw_{-1,\g}(t) = t^{-1}(1-t)^\g$, $0 \le t \le 1$. Then
\begin{align*}
  -  \int_{\VV_0^{d+1}} \Delta_{0,\g} f(x,t) \cdot g(x,t)  \sw_{-1,\g}(x,t)  \d\sm(x,t) 
  = \int_{\VV^{d+1}_0} t(1-t) \frac{\d f}{\d t} \frac{\d g}{\d t} \sw_{-1,\g}(t) \d \sm (x,t) & \\ 
   +   \sum_{1 \le i < j\le d}   \int_{\VV_0^{d+1}} t^{-2}D_{i,j}^{(x)} f(x,t) D_{i,j}^{(x)} g (x,t) \sw_{-1,\g}(t) \d \sm(x,t)&, 
   \end{align*}
where $\frac{\d f}{\d t} = \frac{\d}{\d t} f(t\xi,t)$. In particular, $-\Delta_{0,\g}$ is self-adjoint in $L^2(
\VV^{d+1}, \sw_{-1,\g})$. 
\end{prop} 

The proof of these results will be given in Subsection \ref{sect:proof_corollary}. It is used to deduce the 
following corollary. 

\begin{cor} \label{cor:B-V0_p=2}
Let $d \ge 2$ and $\g > -1$. Then, for $\sw = \sw_{-1,\g}$ and any polynomial $f$,
$$
 \left \| \varphi \partial_t f \right \|_{2,\sw}^2 + \sum_{1\le i< j \le d} \left \|\f{1}{\sqrt{t}} D_{i,j}^{(x)} f \right \|_{2,\sw}^2
       \le  \left \| f \right \|_{2,\sw}  \left \| \Delta_{0,\g} f \right \|_{2,\sw}. 
$$
In particular, for $f \in \Pi_n(\VV_0^{d+1})$, the inequality \eqref{eq:BernsteinLB} with $r =2$ and $p=2$ 
implies the inequalities \eqref{eq:BI-V0-2} and \eqref{eq:BI-V0-3} with $\ell =1$ and $p=2$ when 
$\sw = \sw_{-1,\g}$.  
\end{cor}

In the other direction, the explicit formula of \eqref{eq:LB-operator} and \eqref{eq:D=Dij} shows immediately 
that $\|-\Delta_{0,\g} f\|_{p,\sw}$ is bounded by the sum of $\|\varphi^2 \partial_t^2 f\|_{p,\sw}$, 
$\| \partial_t f\|_{p,\sw}$ and $\| (\frac{1}{\sqrt{t}}D_{i,j})^2 f\|_{p,\sw}$. In particular, the inequalities in the 
Theorem \ref{thm:BI-V0} imply \eqref{eq:BernsteinLB}  with $r =2$. By the spectral property, the inequality 
\eqref{eq:BernsteinLB} is sharp. The corollary and the discussion above provides some assurance that
the inequalities \eqref{eq:BI-V0-2} and \eqref{eq:BI-V0-3} are sharp when $p =2$ and $\sw= \sw_{-1,\g}$. 

\subsection{Main  results on the cone} 
Here the domain is the solid cone in $\RR^{d+1}$ for $d \ge 1$,
$$
\VV^{d+1} = \left \{(x,t): \|x\| \le t, \, 0 \le t \le 1, \, x \in \RR^d \right \}.
$$
Writing $\VV^{d+1}$ as $\{(t \xi, t): \xi \in \sph, 0 \le t \le 1\}$, we see that the first order derivative on the conic 
surfaces are $\partial_t = \f{\partial}{\partial t}$ in the $t$ variable and the angular derivatives $D_{i,j}$, $1\le i, j \le d$, 
in the $x$ variables, which we denote by $D_{i,j}^{(x)}$, as well as one more partial derivative, denoted by
$D_{x_j}$ and defined by 
\begin{equation} \label{eq:Dxj}
     D_{x_j} = \sqrt{t^2-\|x\|^2} \frac{\partial}{\partial x_j}, \qquad 1 \le j \le d, \quad (x,t) \in \VV^{d+1}.
\end{equation}
Let $W$ be a doubling weight on $\VV^{d+1}$; see \ref{sec:OPcone} for the definition. We denote by 
$\|\cdot\|_{p,W}$ the weighted $L^p$ norm on $\VV^{d+1}$
$$
  \|f\|_{p,W} = \left( \int_{\VV^{d+1}} | f(x,t)|^p W(x,t) \d x \d t \right)^{1/p}, \qquad 1 \le p < \infty.  
$$

\begin{thm} \label{thm:BI-V}
Let $W$ be a doubling weight on $\VV^{d+1}$ and let $f\in \Pi_n$. For $\ell \in \NN$ and $1 \le p < \infty$, 
\begin{equation} \label{eq:BI-V-1}
    \left \|\partial_t^\ell  f \right \|_{p,W} \le c_p  n^{2 \ell} \|f\|_{p,W} \quad \hbox{and} \quad
     \left \|\varphi^\ell \partial_t^\ell  f\right \|_{p,W} \le c_p  n^{\ell} \|f\|_{p,W},
 \end{equation}
where $\varphi(t) = \sqrt{t(1-t)}$; moreover, for $ 1 \le  i, j \le d$, 
\begin{equation}\label{eq:BI-V-1B}
   \left \|D_{i,j}^\ell f\right \|_{p,W} \le c_p  n^\ell \|f\|_{p,W},
\end{equation}
and, for $\ell = 1,2$, 
\begin{equation}\label{eq:BI-V-1C}
   \left \|  \frac{1}{\sqrt{t}^\ell } D_{i,j}^\ell f\right \|_{p,W} \le c_p  n^\ell \|f\|_{p,W}, 
\end{equation}
which does not hold for $\ell =3,4,\ldots$ in general. Furthermore, for $\ell \in \NN$ and $1 \le j \le d$, 
\begin{equation}\label{eq:BI-V-2}
   \left \|\partial_{x_j}^\ell  f \right \|_{p,W} \le c_p  n^{2 \ell} \|f\|_{p,W} \quad  \hbox{and} \quad
     \left \|\Phi^\ell \partial_{x_j}^\ell  f\right \|_{p,W} \le c_p  n^{\ell} \|f\|_{p,W},
\end{equation}
where $\Phi(x,t) = \sqrt{t^2-\|x\|^2}$, and, for $\ell = 1, 2$, 
\begin{equation}\label{eq:BI-V-2B}
   \left \|  \frac{1}{\sqrt{t}^\ell } \Phi^\ell \partial_{x_j}^\ell f\right \|_{p,W} \le c_p  n^\ell \|f\|_{p,W}.
\end{equation}
Finally, these inequalities also hold when the norm is the uniform norm $\|\cdot\|_\infty$ on $\VV^{d+1}$. 
\end{thm}
 
In contrast to the inequality \eqref{eq:BI-V-2}, we do not know if \eqref{eq:BI-V-2B} holds for $\ell \ge 3$. As 
it will be shown in Section 4, the two cases are closely related, but the example that works for \eqref{eq:BI-V-2}
does not work for \eqref{eq:BI-V-2B}.
 
Like the case of conic surface, there is also a spectral operator on the cone $\VV^{d+1}$. For $\mu > -\f12$
and $\g> -1$, define the second order differential operator
\begin{align} \label{eq:V-DE}
  \fD_{\mu,\g} : = & \, t(1-t)\partial_t^2 + 2 (1-t) \la x,\nabla_x \ra \partial_t + t \Delta_x -  \langle x, \nabla_x \rangle^2 \\
    &   + (2\mu+d)\partial_t  - (2\mu+\g+d+1)( \la x,\nabla_x\ra + t \partial_t) +   \langle x, \nabla_x \rangle.  \notag
\end{align}
It is proved in \cite{X21} that orthogonal polynomials with respect to the weight function 
$$
W_{\mu,\g}(x,t)= (t^2-\|x\|^2)^{\mu-\f12} (1-t)^\g, \qquad \mu > -\tfrac12, \g > -1
$$ 
on $\VV^{d+1}$ are the eigenfunctions of the operator $\fD_{\mu,\g}$ (see, 
Theorem \ref{thm:Delta0V} below). In particular, the operator satisfies the following Bernstein inequality 
\cite[Theorem 3.1.7]{X21}. 

\begin{thm} \label{thm:BernsteinLB-V}
Let $W$ be a doubling weight on $\VV^{d+1}$. Let $\g > -1$ and $\mu > -\f12$. For $r > 0$ and $1 \le p \le \infty$, 
\begin{equation}\label{eq:BernsteinLB-V}
  \| (-  \fD_{\mu,\g})^{\f r 2}f \|_{p,W} \le c n^r \|f\|_{p,W}, \qquad \forall f\in \Pi_n.
\end{equation}
\end{thm}

The operator $\fD_{\mu,\g}$ is self-adjoint, as can be seen in the following theorem. 

\begin{thm} \label{thm:self-adjointV}
For $\mu > -\f12$ and $\g > -1$, 
\begin{align} \label{eq:intJacobiSolid}
  \int_{\VV^{d+1}}   - \fD_{\mu,\g} f(x,t) & \cdot g(x,t)  W_{\mu,\g}(x,t) \d x \d t =
   \int_{\VV^{d+1}} t \frac{\d f}{\d t} \frac{\d g}{\d t} W_{\mu,\g+1}(x,t) \d x \d t\\
&  + \sum_{i=1}^d \int_{\VV^{d+1}}   D_{x_i} f(x, t) \cdot  
     D_{x_i} g(x,t) t^{-1} W_{\mu,\g}(x) \d x  \d t \notag \\
& + \int_{\VV^{d+1}} \sum_{i < j}   D_{i,j}^{(x)} 
    f(x,t) \cdot D_{i,j}^{(x)} g(x,t) t^{-1} W_{\mu,\g} (x) \d x\,\d t, \notag
\end{align}
where $\frac{\d}{\d t} f = \frac{\d}{\d t} [f(t y,t)]$ for $y \in \BB^d$.
\end{thm} 
 
The identity \eqref{eq:intJacobiSolid} is proved in Subsection \ref{sect:self-adjointV} and it implies immediately 
the following corollary. 
 
\begin{cor}
Let $d \ge 2$, $\mu > -\f12$ and $\g > -1$. Let $W = W_{\mu,\g}$. Then for any polynomial $f$
$$
 \left \| \varphi \partial_t f \right \|_{2,W}^2 + \sum_{j=1}^d \left \|\f{1}{\sqrt{t}} D_{x_j} f\right \|_{2,W}^2 +
     \sum_{1\le i< j \le d} \left \|\f{1}{\sqrt{t}} D_{i,j}^{(x)} f \right \|_{2,W}^2
       \le  \left \| f \right \|_{2,W}  \left \| \fD_{\mu,\g} f \right \|_{2,W}. 
$$
In particular, the inequality \eqref{eq:BernsteinLB-V} for $\fD_{\mu,\g}$ with $r=2$ implies the inequalities 
\eqref{eq:BI-V-2} for $\varphi \partial_t$ and \eqref{eq:BI-V0-3} for $D_{i,j}^{(x)}$ and $D_{x_j}$ with 
$\ell =1$ and $p=2$ for $W = W_{\mu,\g}$. 
\end{cor}

The corollary follows as in Corollary \ref{cor:B-V0_p=2}. The Bernstein inequalities for the first order derivatives 
in Theorem \ref{thm:BI-V} imply \eqref{eq:BernsteinLB-V} when $r =2$. Together, they provide some assurance 
that the inequalities in Theorem \ref{thm:BI-V} are sharp. 

\subsection{Bernstein inequality on the triangle}
For $d =1$, the cone becomes the triangle $\VV^2 = \{(x,t): 0 \le t \le 1, \, |x| \le t\}$ of $\RR^2$. Making an
affine change of variable $(x,t) = (y_1-y_2, y_1 + y_2)$, the triangle $\VV^2$ becomes the standard triangle
domain
$$
     \TT^2 = \{(y_1, y_2): y_1 \ge 0, \, y_2 \ge 0, \,  y_1+ y_2 \le 1\}.
$$
The triangle $\TT^2$ is symmetric under permutation of $\{y_1,y_2,1-y_1-y_2\}$ and it is customary to consider
the derivatives (cf. \cite{BX, DX}) 
$$
  \partial_1 = \partial_{y_1}, \quad \partial_2 = \partial_{y_2}, \quad \partial_3 = \partial_{y_2} - \partial_{y_1}.
$$
We further define, for $y = (y_1,y_2) \in \TT^2$,  
$$
\phi_1(y) = \sqrt{y_1(1-y_1-y_2)}, \quad \phi_2(y) = \sqrt{y_2(1-y_1-y_2)}, \quad \phi_3(y) = \sqrt{y_1 y_2}. 
$$
Let $\|\cdot\|_{p,w}$ be the $L^p$ norm on the triangle $\TT^2$ for $1\le p < \infty$ and the uniform norm on $\TT^2$
if $p = \infty$. 

\begin{thm} 
Let $w$ be a doubling weigh on $\TT^2$. Let $f \in \Pi_n^2$. Then, for $1 \le p < \infty$, $\ell \in \NN$ and
$i =1,2,3$, 
\begin{equation} \label{eq:triangle}
 \left \| \partial_i^\ell f\right\| _{p,w} \le c n^{2\ell} \|f \|_{p,w}\quad  \hbox{and}\quad 
  \left \| \phi_i^\ell \partial_i^\ell f\right\| _{p,w} \le c n^\ell \|f\|_{p,w}.
\end{equation} 
Furthermore, for $\ell =1,2$, 
\begin{align} 
 \left \| \frac{1}{(1-y_2)^{\ell/2}} \phi_1^\ell  \partial_{1}^\ell f \right  \|_{p,w} \le c n^\ell \|f\|_{p,w},   \label{eq:tri1}\\  
 \left \| \frac{1}{(1-y_1)^{\ell/2}} \phi_2^\ell  \partial_{2}^\ell f \right  \|_{p,w} \le c n^\ell \|f\|_{p,w},    \label{eq:tri2}\\ 
 \left \| \frac{1}{(y_1+y_2)^{\ell/2}} \phi_3^\ell  \partial_{3}^\ell f \right  \|_{p,w} \le c n^\ell \|f\|_{p,w}. \label{eq:tri3}
\end{align}
Moreover, these inequalities hold when the norm is replaced by the uniform norm. 
\end{thm}
 
\begin{proof}
Under the change of variables $\VV^2 \mapsto \TT^2$: $(x,t) \mapsto (y_1-y_2, y_1 + y_2)$, we see that 
$$
\partial_x= \frac12 (\partial_{y_2} - \partial_{y_1}) = \frac12 \partial_3 \quad\hbox{and} \quad 
\sqrt{t^2-|x|^2} = 2 \sqrt{y_1 y_2} = 2 \phi_3(y). 
$$
Hence, for $i = 3$, the two inequalities in \eqref{eq:triangle} follow from \eqref{eq:BI-V-2}, where $d=1$, 
and the inequality \eqref{eq:tri3} follows from \eqref{eq:BI-V-2B} in Theorem \ref{thm:BI-V}. The inequalities 
for $\partial_1$ and $\partial_2$ follow from those
for $\partial_3$ by a change of variables that amounts to permute the variables $(y_1,y_2,1-y_1-y_2)$. 
\end{proof}

The inequality \eqref{eq:triangle} is known when $w$ is a constant weight (cf. \cite{BX, DT, T1}) or the Jacobi 
weight. It is accepted as a natural generalization of the Bernstein inequality \eqref{eq:Bernstein[0,1]phi} 
of one variable, since the left-hand side of the inequality in the uniform norm gives a pointwise inequality 
that reduces to \eqref{eq:Bernstein[0,1]phi} when restricted to the boundary of the triangle. The latter property, 
however, is satisfied by \eqref{eq:tri1}--\eqref{eq:tri3} as well. Moreover, the factors in the left-hand 
side remains bounded; for example, $1-y_1-y_2 \le 1 - y_2$ in \eqref{eq:tri1}. The inequalities 
\eqref{eq:tri1}--\eqref{eq:tri3} are surprising and they are new as far as we are aware. It is also suggestive
and may be worthwhile to ask if the same phenomenon appears on other domains, such as polytopes
\cite{T1,T2}. 

\section{Bernstein inequalities on conic surface}\label{sec:OP-conicS}
\setcounter{equation}{0}

We work on the conic surface in this section. In the first subsection, we recall what is needed for our analysis
on the domain, based on orthogonal polynomials and highly localized kernels. Several auxiliary inequalities
are recalled or proved in the second subsection, to be used in the rest of the section. The proofs of the main 
results are based on the estimates of the kernel functions, which are carried out in the third subsection
for the derivatives in the $t$-variable and the fourth subsection for the angular derivatives. Finally, the proof 
of Theorem \ref{thm:self-adjointV} and its corollary is given in the fifth subsection.

\subsection{Orthogonal polynomials and localized kernels}\label{sec:OPconicSurface}
Let $\Pi(\VV_0^{d+1})$ denote the space of polynomials restricted on the conic surface $\VV_0^{d+1}$
and, for $n = 0,1,2,\ldots$, let $\Pi_n(\VV_0^{d+1})$ be the subspace of polynomials in $\Pi(\VV_0^{d+1})$
of total degree at most $n$. Since $\VV_0^{d+1}$ is a quadratic surface, it is known that 
$$
    \dim \Pi_n(\VV_0^{d+1}) = \binom{n+d}{n}+\binom{n+d-1}{n-1}.
$$
For $\b > - d$ and $\g > -1$, we define the weight function $\sw_{\b,\g}$ by
$$
   \sw_{\b,\g}(t) = t^\b (1-t)^\g, \qquad 0 \le t \le 1.
$$
Orthogonal polynomials with respect to $ \sw_{\b,\g}$ on $\VV_0^{d+1}$ are studied in \cite{X20}. Let 
$$
\la f, g\ra_{\sw_{\b,\g}} =\bs_{\b,\g} \int_{\VV_0^{d+1}} f(x,t) g(x,t) \sw_{\b,\g} \d \sm(x,t),
$$ 
where $\d \sm$ denotes the Lebesgue measure on the conic surface, which is a well defined inner product 
on $\Pi(\VV_0^{d+1})$. Let $\CV_n(\VV_0^{d+1},\sw_{\b,\g})$ be the space of orthogonal polynomials of 
degree $n$. Then $\dim \CV_0(\VV_0^{d+1},\sw_{\b,\g}) =1$ and 
$$
   \dim \CV_n(\VV_0^{d+1},\sw_{\b,\g})  = \binom{n+d-1}{n}+\binom{n+d-2}{n-1},\quad n=1,2,3,\ldots.
$$
Let $\CH_m(\sph)$ be the space of spherical harmonics of degree $m$ in $d$ variables. Let 
$\{Y_\ell^m: 1 \le \ell \le \dim \CH_m(\sph)\}$ denote an orthonormal basis of $\CH_m(\sph)$. 
Then the polynomials
\begin{equation*} 
  \sS_{m, \ell}^n (x,t) = P_{n-m}^{(2m + \b + d-1,\g)} (1-2t) Y_\ell^m (x), \quad 0 \le m \le n, \,\, 
      1 \le \ell \le \dim \CH_m(\sph),
\end{equation*}
consist of an orthogonal basis of $\CV_n(\VV_0^{d+1}, \sw_{\b,\g})$. Let $\Delta_{0,\g}$ be the differential
operator defined in \eqref{eq:LB-operator}. It has orthogonal polynomials as eigenfunctions \cite[Theorem 7.2]{X20}.

\begin{thm}\label{thm:Jacobi-DE-V0}
Let $d\ge 2$ and $\g > -1$. The orthogonal polynomials in $\CV_n(\VV_0^{d+1}, \sw_{-1,\g})$ are eigenfunctions
of $\Delta_{0,\g}$; more precisely, 
\begin{equation}\label{eq:eigen-eqn}
    \Delta_{0,\g} u =  -n (n+\g+d-1) u, \qquad \forall u \in \CV_n(\VV_0^{d+1}, \sw_{-1,\g}).
\end{equation}
\end{thm} 

The reproducing kernel of the space $\CV_n(\VV_0^{d+1}, \sw_{\b,\g})$ is denoted by $\sP_n(\sw_{\b,\g};\cdot,\cdot)$, 
which can be written as
$$
\sP_n\big(\sw_{\b,\g}; (x,t),(y,s) \big) = \sum_{m=0}^n \sum_{k=1}^{\dim \CH_m^d}
    \frac{  \sS_{m, \ell}^n(x,t)  \sS_{m, \ell}^n(y,s)}{\la  \sS_{m, \ell}^n,  \sS_{m', \ell'}^{n'} \ra_{\sw_{\b,\g}}}. 
$$
Let $\proj_n(\sw_{\b,\g}): L^2(\VV_0^{d+1},\sw_{\b,\g}) \to \CV_n(\VV_0^{d+1}, \sw_{\b,\g})$ be the orthogonal 
projection operator. It is an integral operator with $\sP_n\big(\sw_{\b,\g}; \cdot,\cdot \big)$ as its kernel, 
$$
\proj_n(\sw_{\b,\g};f) = \int_{\VV_0^{d+1}} f(y,s) \sP_n\big(\sw_{\b,\g}; \,\cdot, (y,s) \big)  \sw_{\b,\g}(s) \d\sm(y,s).
$$
The kernel $ \sP_n\big(\sw_{\b,\g}; \,\cdot, \cdot \big)$ satisfies a closed form formula, called the addition formula 
since it is akin to the classical addition formula for the spherical harmonics. The addition formula is given in terms 
of the Jacobi polynomial $P_n^{(\a,\b)}$, the orthogonal polynomial with respect to $w_{\a,\b}(t) = (1-t)^\a(1+t)^\b$ 
on $[-1,1]$. Let 
$$
  Z_{n}^{(\a,\b)}(t) = \frac{P_n^{(a,\b)}(t) P_n^{(\a,\b)}(1)}{h_n^{(\a,\b)}},  \quad \a,\b > -1,
$$
where $h_n^{(\a,\b)}$ is the $L^2$ norm of $P_n^{(\a,\b)}$ in $L^2([-1,1],w_{\a,\b})$. The addition formula is
of the simplest form when $\g = -1$. 

\begin{thm}  \label{thm:sfPbCone2}
Let $d \ge 2$ and $\g \ge -\f12$. Then, for $(x,t), (y,s) \in \VV_0^{d+1}$,
\begin{align} \label{eq:sfPbCone}
 \sP_n \big(\sw_{-1,\g}; (x,t), (y,s)\big) =  b_{\g,d}  \int_{[-1,1]^2} & Z_{n}^{(\g+d-\f32,-\f12)} \big( 2 \zeta (x,t,y,s; v)^2-1\big) \\
  & \times  (1-v_1^2)^{\f{d-4}{2}} (1-v_2^2)^{\g-\f12} \d v, \notag
\end{align} 
where $b_{\g,d}$ is a constant so that $\sP_0\big(\sw_{-1,\g}; (x,t), (y,s)\big) =1$ and 
\begin{equation}\label{eq:zetaV0}
 \zeta (x,t,y,s; v)  = v_1 \sqrt{\tfrac{st + \la x,y \ra}2}+ v_2 \sqrt{1-t}\sqrt{1-s};
\end{equation}
moreover, the identity holds under limit when $\g = -\f12$ and/or $d = 2$. 
\end{thm} 

Our main tool for establishing polynomial inequalities is the highly localized kernel defined via a smooth
cur-off function $\wh a \in C^\infty(\RR)$, which is a non-negative function and satisfies 
$\mathrm{supp}\, \wh a \subset [0, 2]$ and $\wh a(t) = 1$, $t\in [0, 1]$. The kernel is defined by 
\begin{equation} \label{def:Ln-gen}
   \sL_n (\sw_{-1,\g}; (x,t),(y,s)) = \sum_{k=0}^{\infty} \wh a \left(\frac{k}{n}\right) \sP_n(\sw; (x,t),(y,s)). 
\end{equation}
Since $\wh a$ is supported on $[0,2]$, this is a kernel of polynomials of degree at most $2n$ in either
the $x$ or the $y$ variable. Using the closed form of the reproducing kernel, we can write $\sL(\sw_{\b,\g})$ 
in terms of the kernel of the Jacobi polynomials defined by
$$
   L_n^{(\l,-\f12)} (t) = \sum_{k=0}^\infty\wh a \left(\frac{k}{n}\right) Z_n^{(\l,-\f12)}(t). 
$$
Indeed, it follows immediately from \eqref{eq:sfPbCone} that 
\begin{align}\label{eq:Ln-intV0}
\sL_n (\sw_{-1,\g}; (x,t), (y,s) )=  
    b_{\g,d}  \int_{[-1,1]^2} & L_n ^{(\g+d-\f32,-\f12)}\big(2 \zeta (x,t,y,s; v)^2-1 \big)\\
  &  \times    (1-v_1^2)^{\f{d-2}2-1}(1-v_2^2)^{\g-\f12} \d v. \notag
\end{align}

To show that this kernel is highly localized, we need the distance on the conic surface. For $(x,t)$ 
and $(y,s)$ on $\VV_0^{d+1}$, the distance function $\sd_{\VV_0}$ on $\VV_0^{d+1}$ is defined by 
\begin{equation}\label{eq:distV0}
  \sd_{\VV_0} ((x,t), (y,s)): =  \arccos \left(\sqrt{\frac{\la x,y\ra + t s}{2}} + \sqrt{1-t}\sqrt{1-s}\right).
\end{equation}
For $r > 0$ and $(x,t)$ on $\VV_0^{d+1}$, we let $\sc((x,t), r)$ be the ball centered at $(x,t)$ with radius $r$
in terms of this distance function; that is, 
$$
      \sc((x,t), r): = \left\{ (y,s) \in \VV_0^{d+1}: \sd_{\VV_0} \big((x,t),(y,s)\big)\le r \right\}.
$$   
Let $E$ be a subset in $\VV_0^{d+1}$ and $\sm$ denotes the Lebesgue measure in $\VV_0^{d+1}$. We define 
$$
  \sw (E) = \int_E \sw(x,t) \d \sm(x,t).
$$
A weight function $\sw$ is a doubling weight if there is a constant $L > 0$ such that 
$$
   \sw\big(\sc((x,t), 2 r)\big) \le L \, \sw\big(\sc((x,t), r)\big), \quad r >0.
$$
The least constant $L$ is called a doubling constant and the doubling index $\a(\sw)$ is the least index for 
which $\sup_{\sc((x,t), r} \sw(\sc((x,t), 2^m r))) /\sw(\sc((x,t), r)) \le c_{L(\sw)} 2^{m \a (\sw)}$, $m=1,2,\ldots$. 
As an example, the weight $\sw_{\b,\g}$ is a doubling weight on $\VV_0^{d+1}$ \cite[Proposition 4.6]{X21}
and, for $n = 1,2, \ldots$, $\sw_{-1,\g}\big(\sc((x,t), n^{-2})\big)\sim \sw_{\g,d} (n; t)$ with
\begin{equation}\label{eq:w(n;t)}
     \sw_{\g,d} (n; t) = \big(1-t+n^{-2}\big)^{\g+\f12}\big(t+n^{-2}\big)^{\f{d-2}{2}}. 
\end{equation}

It is shown in \cite[Theorem 4.10]{X21} that, for $d\ge 2$, $\g \ge -\f12$, and any $\k > 0$, the 
kernel $\sL_n (\sw_{-1,\g}; \cdot,\cdot)$ satisfies the estimate 
\begin{equation}\label{eq:L-localized}
\left |\sL_n (\sw_{-1,\g}; (x,t), (y,s))\right|
\le \frac{c_\k n^d}{\sqrt{ \sw_{\g,d} (n; t) }\sqrt{ \sw_{\g,d} (n; s) }}
\big(1 + n \sd_{\VV_0}( (x,t), (y,s)) \big)^{-\k}.
\end{equation}
This shows, in particular, that the kernel decays away from $(x,t) = (y,s)$ faster than any polynomial rate. The 
decaying estimate is established using the case of $m = 0$ in the following lemma  \cite[Theorem 2.6.7]{DaiX}. 

\begin{lem}
Let $\ell$ be a positive integer and let $\eta$ be a function that satisfy, $\eta\in C^{3\ell-1}(\RR)$, 
$\mathrm{supp}\, \eta \subset [0,2]$ and $\eta^{(j)} (0) = 0$ for $j = 0,1,2,\ldots, 3 \ell-2$. Then, 
for $\a \ge \b \ge -\f12$, $t \in [-1,1]$ and $n\in \NN$, 
\begin{equation} \label{eq:DLn(t,1)}
\left| \frac{d^m}{dt^m} L_n^{(\a,\b)}(t) \right|  \le c_{\ell,m,\a}\left\|\eta^{(3\ell-1)}\right\|_\infty 
    \frac{n^{2 \a + 2m+2}}{(1+n\sqrt{1-t})^{\ell}}, \quad m=0,1,2,\ldots. 
\end{equation}
\end{lem}

The proof of the localization of the kernel requires further properties of the distance function. 
Let $\sd_{[-1,1]}(\cdot,\cdot)$ be the distance function on the interval $[0,1]$ defined by
$$
\sd_{[-1,1]}(t, s)  = \arccos \left(\sqrt{t}\sqrt{s} + \sqrt{1-t}\sqrt{1-s}\right), \qquad t, s \in [0,1]
$$
and let $\sd_{\SS}(\cdot,\cdot)$ denote the geodesic distance on the unit sphere $\sph$ defined by  
$$
  \sd_\SS (\xi,\eta) = \arccos \la \xi,\eta\ra, \qquad \xi, \eta \in \sph.
$$
Then, for $(x,t), (y,s) \in \VV_0^{d+1}$ and setting $x= t \xi$ and $y = s \eta$, the distance $\sd_{\VV_0}(\cdot,\cdot)$
satisfies, as shown in \cite[Proposition 4.3]{X21},
\begin{equation} \label{eq:d2=d2+d2}
   c_1 \sd_{\VV_0} ((x,t), (y,s)) \le \sd_{[0,1]}(t,s)  + (t s )^{\f14} \sd_{\SS}(\xi,\eta) \le  c_2 \sd_{\VV_0} ((x,t), (y,s)).
\end{equation}
Another two inequalities that we shall need are
\begin{equation} \label{eq:|s-t|}
  \big| \sqrt{t} - \sqrt{s} \big|\le \sd_{\VV_0} ((x,t), (y,s)) \quad \hbox{and} \quad 
      \big| \sqrt{1-t} - \sqrt{1-s} \big| \le \sd_{\VV_0} ((x,t), (y,s))
\end{equation}
for $(x,t), (y,s) \in \VV_0^{d+1}$, established in \cite[Lemma 4.4]{X21}.

\subsection{Auxiliary inequalities} 

We state several auxiliary results in this subsection. First of all, we will need the maximal function $f_{\b,n}^\ast$ 
defined by
\begin{equation} \label{eq:fbn*}
f_{\b,n}^\ast(x,t) = \max_{(y,s)\in \VV_0^{d+1}} \frac{|f(y,s)|} {\left(1+n \sd_{\VV_0}((x,t),(y,s)) \right)^\b}.
\end{equation}
This maximal function satisfies the following property \cite[Corollary 2.11]{X21}.

\begin{prop}\label{prop:fbn-bound}
If $ 0< p\leq \infty$, $ f\in\Pi_n(\VV_0^{d+1})$ and $\b > \a(\sw)/p$, then
\begin{equation} \label{eq:fbn*bound}
  \|f\|_{p, \sw} \leq \|f_{\b,n}^\ast\|_{p,\sw} \leq c  \|f\|_{p,\sw},
\end{equation}
where $c$ depends also on $L(\sw)$ and $\b$ when $\b$ is either large or close to $\a(\sw)/p$.
\end{prop}

For our next result, we need the concept of a maximal $\ve$-separated set on the conic surface. However,
the concept will be needed for the unit ball and the solid cone later in the paper, so we give its definition
on a domain $\Omega$ equipped with a distance function $\sd(\cdot,\cdot)$. Let 
$B(x,\ve) = \{y \in \Omega:\sd(x,y) \le \ve\}$ be the ball centered at $x$ and with radius $\ve$ in $\Omega$. 

\begin{defn}\label{defn:separated-pts}
Let $\Xi$ be a discrete set in $\Omega$. 
\begin{enumerate} [  \quad (a)]
\item Let $\ve>0$. A discrete subset $\Xi$ of $\Omega$ is called $\ve$-separated if $\sd(x,y) \ge\ve$
for every two distinct points $x, y \in \Xi$. 
\item $\Xi$ is called maximal if there is a constant $c_d > 1$ such that 
\begin{equation*} 
  1 \le  \sum_{z\in \Xi} \chi_{B(z, \ve)}(x) \le c_d, \qquad \forall x \in \Omega,
\end{equation*}
where $\chi_E$ denotes the characteristic function of the set $E$.
\end{enumerate}
\end{defn} 

A maximal $\ve$-separated subset is constructed in \cite[Proposition 4.17]{X21} using the separated sets 
in the $t$-variable on $[0,1]$ and in the $\xi$ variable on $\sph$, where $(x,t) = (t \xi, t) \in \VV_0^{d+1}$. We 
recall what is necessary for our purpose. Let $\ve > 0$ and let $N = \lfloor \frac{\pi}{2}\ve^{-1} \rfloor$. We
define
\begin{equation} \label{eq:tj-epj}
  t_j = \sin^2  \frac{(2j-1)\pi}{4 N} \quad \hbox{and}\quad \ve_j =  \frac{\pi \ve} {2 \sqrt{t_j}}, \quad 1 \le j \le N.
\end{equation}
Let $\Xi_\SS(\ve_j)$ be the maximal $\ve_j$-separated set of $\sph$, so that there is a family of sets 
$\{\SS_\xi(\ve_j): \xi \in \Xi_\SS(\ve_j)\}$ which forms a partition $\sph = \bigcup_{\eta \in \Xi_\SS(\ve_j)} \SS_\eta(\ve_j)$. 
Then 
\begin{equation} \label{eq:separateV0}
   \Xi_{\VV_0} = \big\{(t_j \xi, t_j): \,  \xi \in \Xi_\SS(\ve_j), \, 1\le j \le N \big\} 
\end{equation}
defined a maximal $\ve$-separated subset on $\VV_0^{d+1}$. Such a set is used to establish the
Marcinkiewicz-Zygmund inequality on the conic surface, which we need below. 

Out next result is of interest in itself. For polynomials of one variable, it is known \cite[(7.1.7)]{MT1}  
$$
   \int_{-1}^1 |f(t)|^p w(t) \d t \le c_\delta \int_{-1+ \delta n^{-2}}^{1-\delta n^{-2}} |f(t)|^p w(t) \d t, 
   \quad \deg f \le n,
$$
where $w$ is a doubling weight and $\delta$ is a positive constant. The following proposition is an 
analog of the above inequality on the conic surface. 

\begin{prop} \label{prop:Remz}
Let $\sw$ be a doubling weight function on $\VV_0^{d+1}$. For $n \in \NN$, let $\chi_{n,\delta}(t)$ denote the 
characteristic function of the interval $[\f \delta{n^2}, 1- \f \delta {n^2}]$. Then, for $f \in \Pi_n(\VV_0^{d+1})$,
$1 \le p < \infty$, and every $\delta > 0$, 
\begin{equation}\label{eq:Remz}
  \int_{\VV_0^{d+1}} |f(x,t)|^p \sw(x,t) \d \sm(x,t) \le c_\delta \int_{\VV_{0}^{d+1}} |f(x,t)|^p \chi_{n,\delta}(t) \sw(x,t) \d \sm(x,t). 
\end{equation}
Moreover, when $p = \infty$, 
\begin{equation}\label{eq:Remz2}
  \|f\|_\infty \le c \|f \chi_{n,\delta}(t)\|_{\infty}.
\end{equation}
\end{prop}

\begin{proof}
The proof uses the Marcinkiewicz-Zygmund inequality on the conic surface. Let $\Xi_{\VV_0}$ be
the $\ve$-maximal separated set in \eqref{eq:separateV0}. For $\ve = \b n^{-1}$ with $\b > 0$, the 
Marcinkiewicz-Zygmund inequality states \cite[Theorem 4.18]{X21} that, for $1 \le p < \infty$, 
\begin{enumerate}[$(i)$]
\item for $f\in\Pi_m(\VV_0^{d+1})$ with $n \le m \le c n$,
\begin{equation*}
  \sum_{(z,r) \in \Xi_{\VV_0}} \Big( \max_{(x,t)\in \sc \big((z,r), n^{-1}\big)} |f(x,t)|^p \Big)
     \sw\!\left(\sc \big((z, r), n^{-1}\big) \right) \leq c_{\sw} \|f\|_{p,\sw}^p;
\end{equation*}
\item for 
$f \in \Pi_n(\VV_0^{d+1})$,  
\begin{align*}
  \|f\|_{p,\sw}^p \le c_{\sw} \sum_{(z,r) \in\Xi_{\VV_0}} \Big(\min_{(x,t)\in \sc\bigl((z,r), n^{-1}\bigr)} |f(x,t)|^p\Big)
          \sw\bigl(\sc\big((z,r),n^{-1}\big)\bigr).
\end{align*}
\end{enumerate}
Clearly $t_0 = \sin^2 \frac{\pi}{4 N} \sim n^{-2}$ and $1- t_N \sim n^{-2}$ and the constant is proportional to 
$\b^{-1}$. Hence, for a fixed $\delta > 0$, by choosing $\b$ sufficiently large, we see that $f(z,r) =
 f(z,r) \chi_{n,\delta}(r)$ for all $(z,r) \in \Xi_{\VV_0}$. Consequently, it follows from (ii) and (i) that
\begin{align*}
     \|f\|_{p,\sw}^p \, & \le c_{\sw} \sum_{(z,r) \in\Xi_{\VV_0}} |f(z,r) \chi_{n,\delta}(r) |^p
          \sw\bigl(\sc\big((z,r),n^{-1}\big)\bigr) \\
            & \le c \int_{\VV_{0}^{d+1}} |f(x,t)|^p \chi_{n,\delta}(t) \sw(x,t) \d \sm(x,t)
\end{align*}
for all $f\in \Pi_n(\VV_0^{d+1})$, which is the desired inequality \eqref{eq:Remz}. For the uniform norm, 
$$
  |f(x,t)| \le c \max_{(z,r) \in\Xi_{\VV_0}} |f(z,r)|
$$
by \cite[Theorem 2.14]{X21}, from which the inequality \eqref{eq:Remz2} follows readily. 
\end{proof}

\begin{rem}\label{rem:remark1}
Setting $x = t\xi$, $\xi \in \sph$ shows that the space $\Pi_n^d(\VV_0^{d+1})$ is a subspace of $\Pi_n^*$, which
contains all polynomials of degree $n$ in $t$ and in $\xi$. Examining the proof in \cite{X21} carefully, it is not 
difficult to see that the Marcinkiewicz-Zygmund inequalities in \cite[Theorem 4.18]{X21} holds for $\Pi_n^*$ and, 
as a consequence, so does \eqref{eq:Remz}. 
\end{rem}

Our next lemma is technical and reduces to \cite[Lemma 4.14]{X21} when $\rho =0$
and $\tau_1=\tau_2=0$. The general version is needed to handle the derivatives of the highly localized kernels. 

\begin{lem}\label{lem:intLn}
Let $d\ge 2$ and $\g > -1$. Let $\rho$ satisfies $0 \le \rho \le (d-2)/2$, $\tau_1$ satisfies 
$0 \le \tau_1 \le 1/2$ and $\tau_2$ satisfies $0 \le \tau_2 \le d/2$. For $0 < p < \infty$, 
assume $\k > \frac{2d}{p} + (\g+\f{d-1}{2}) |\f1p-\f12|$. Then for $x = t \xi$ and $y = s \eta$ and 
$(x,t) \in \VV_0^{d+1}$,  
\begin{align}\label{eq:intLn1}
\int_{\VV_0^{d+1}} \frac{ s^{-\tau_1}(1-s)^{-\tau_2}
     (1 + \la \xi, \eta \ra)^{-\rho} \sw_{-1,\g}(s)  \d \sm(y,s) }{  \sw_{\g,d} (n; s)^{\f{p}2}
    \big(1 + n \sd_{\VV_0}( (x,t), (y,s)) \big)^{\k p}} 
    \le c n^{-d +2 |\tau|} \sw_{\g,d} (n; t)^{1-\f{p}{2}}.
\end{align}
\end{lem}

\begin{proof}
For $\rho =0$ and $\tau_i =0$, this estimated is established in \cite[Lemma 4,1,4]{X21}. The proof follows 
along the same argument, we provide an account on the modification. The first step is to show that it is 
sufficient to consider only $p=2$ by using the doubling property of $\sw_{-1,\g}$, which remains valid when 
$\rho > 0$ and $\tau_i > 0$. Let $J_2$ denote the left-hand side of \eqref{eq:intLn1} 
with $p=2$. Then, using 
\begin{equation} \label{eq:intSS}
  \int_{\sph} g(\la \xi,\eta\ra) d\s(\eta) 
      = \o_{d-1} \int_{-1}^1 g(u) (1-u^2)^{\f{d-3}{2}} \d u
\end{equation} 
and $\t \sim 2 \sin\f{\t}2 = \sqrt{1-\cos \t}$ for $\t \le [0,\pi]$, we deduce that 
\begin{align*}
J_{2} \, &\le c  \int_0^1 \int_{-1}^1 \frac{ s^{d-1-\tau_1-\rho} (1-s)^{-\tau_2} (1+u)^{-\rho}(1-u^2)^{\f{d-3}{2}}\sw_{-1,\g}(s)}
{\sw_{\g,d} (n; s)  \big(1 + n \arccos \big( \sqrt{ts} \sqrt{\frac{1+u}{2}} + \sqrt{1-t}\sqrt{1-s} \big)\big)^{2\k} } \d u  \d s \\
   & \le c  \int_0^1 \int_{0}^1 \frac{s^{d-1-\tau_1-\rho} (1-s)^{-\tau_2} v^{d-2-2\rho}(1-v^2)^{\f{d-3}{2}}\sw_{-1,\g}(s) }
          {\sw_{\g,d} (n; s) \left(1 + n \sqrt{1- \sqrt{ts} v -  \sqrt{1-t}\sqrt{1-s}} \right)^{2k}} \d v  \d s,
\end{align*}
where the second step follows from the changing variable $ \sqrt{\frac{1+u}{2}} \mapsto v$. Making a further changing 
of variable $v \mapsto z/\sqrt{s}$ gives 
\begin{align*}
  J_{2} \,& \le c    \int_0^1 \int_{0}^{\sqrt{s}} \frac{s^{1-\tau_1+\rho} (1-s)^{-\tau_2}\, z^{d-2-2\rho}(s- z^2)^{\f{d-3}{2}}
      \sw_{-1,\g}(s) }{\sw_{\g,d} (n; s) \left(1 + n \sqrt{1- \sqrt{t}\,z -  \sqrt{1-t}\sqrt{1-s}} \right)^{2\k} }  \d z  \d s.
\end{align*}       
We now use $s^{1-\tau_1+\rho} z^{d-2 -2 \rho} \le s^{\frac{d}{2} - \tau_1} \le (s+n^{-2})^{\frac{d}{2}-\tau_1} 
\le n^{2 \tau_1}(s+n^{-2})^{\frac{d}{2}}$, which holds for $\rho \le (d-2)/2$ and $\tau_1 \le d/2$, and the 
formulas for $\sw_{-1,\g}(s)$ and $\sw_{\g,d} (n; s)$ to deduce  
\begin{align*}
J_2  \le c n^{2 \tau_1}  \int_0^1 \int_{0}^{\sqrt{s}} \frac{ (1-s)^{-\tau_2}(s- z^2)^{\f{d-3}{2}} }{ (1-s+n^{-2})^{\f12} 
       \left(1 + n \sqrt{1- \sqrt{t}\, z -  \sqrt{1-t}\sqrt{1-s}} \right)^{2\k} } \d z  \d s.        
\end{align*}
We make one more change of variable $s\mapsto 1-w^2$ and make use of $1-s = w^2/ (1+\sqrt{1-w^2}) \le w^2$ to obtain
\begin{align*}
   J_{2}  \, & \le c n^{2 \tau_1}    \int_0^1 \int_{0}^{\sqrt{1-w^2}}  \frac{w^{1-2\tau_2} (1-w^2- z^2)^{\f{d-3}{2}} }{
   \left(w^2+n^{-2}\right)^{\f12} \left(1 + n \sqrt{1- \sqrt{t} z -  \sqrt{1-t} \, w} \right)^{2\k}} \d z  \d w \\
  \, &  \le c n^{2 \tau_1+2\tau_2}   \int_0^1 \int_{0}^{\sqrt{1-w^2}}  \frac{(1-w^2- z^2)^{\f{d-3}{2}} }
        {\left(1 + n \sqrt{1- \sqrt{t} z -  \sqrt{1-t} \, w} \right)^{2\k}} \d z  \d w, 
\end{align*}
where the second inequality follows from $w^{1-2\tau_2} \le (w^2+n^{-2})^{\f12 - \tau_2} \le n^{2\tau_2}(w^2+n^{-2})^{\f12}$,
which holds for $\tau \le \f12$. The last integral in the right-hand side has appeared and estimated in the proof 
of \cite[Lemma 4.14]{X21}, from which the desired estimate follows. 
\end{proof}

Finally we recall the following lemma \cite[Lemma 4.11]{X21}, which plays an important role in establishing the localization
of the kernel $\sL_n(\sw_{-1,\g})$. 

\begin{lem} \label{lem:kernelV0}
Let $d \ge 2$ and $\g > -\f12$. Then, for $\b \ge 2\g+ d+1$, 
\begin{align*}
  \int_{[-1,1]^2} & \frac{(1-v_1)^{\f{d-2}2-1}(1-v_2)^{\g-\f12}} 
        {\big(1+n\sqrt{1- \zeta (x,t,y,s; v)}\,\big)^{\b}} \d v\\
    & \qquad \le \frac{cn^{- (2\g+d-1)}} 
   {\sqrt{\sw_{\g,d}(n; t)}\sqrt{\sw_{\g,d}(n; s)}\big(1+n \sd_{\VV_0}((x,t),(y,s))\big)^{\b - 3\g- \frac{3d+1}{2}}} .
\end{align*}
\end{lem}

\subsection{Bernstein inequality for $\partial_t$ on the conic surface}\label{sec:BernsteinV0}

To prove the Bernstein inequality for the first-order derivatives, we need to understand the action
of these derivatives on the highly localized kernels. We start with the case $\partial_t$. The key ingredient 
lies in the estimates below. 

\begin{lem} \label{lem:D2kernelV0}
Let $d\ge 2$ and $\g \ge -\f12$. Then for any $\k > 0$,
 \begin{equation*} 
\left | \partial_t \sL_n (\sw_{-1,\g}; (x,t), (y,s))\right|
\le \frac{c_\k n^{d+1}}{ \varphi(t) \sqrt{ \sw_{\g,d} (n; t) }\sqrt{ \sw_{\g,d} (n; s) }
\big(1 + n \sd_{\VV_0}( (x,t), (y,s)) \big)^{\k}}.  
\end{equation*}
\end{lem}

\begin{proof}
Taking derivative on the integral expression \eqref{eq:Ln-intV0} of the kernel and writing $\zeta(v)$ in place of 
$\zeta (x,t,y,s; v)$, we obtain with $\l = \g+ d-1$,  
\begin{align}\label{eq:DLn-intV0}
   \partial_t \sL_n (\sw_{-1,\g}; (x,t), (y,s) ) & =  
    c_{\g}  \int_{[-1,1]^2}  \partial L_n ^{(\l-\f12,-\f12)}\big(2 \zeta (v)^2-1 \big)\\
  & 
  \times 4  \zeta (v) \partial_t \zeta (v) (1-v_1^2)^{\f{d-2}2-1}(1-v_2^2)^{\g-\f12} \d v, \notag
\end{align}
where $\partial L_n ^{(\l-\f12,-\f12)}(z) = \f{\sd}{\sd z} L_n ^{(\l-\f12,-\f12)}(z)$. We first prove the following estimate
\begin{align}\label{eq:Dt-zeta}
  \left|2 \varphi(t) \partial_t \zeta (x,t,y,s; v)\right| \le \Sigma_1 + \Sigma_2(v_1) + \Sigma_3(v_2),
\end{align}  
where 
$$  
   \Sigma_1= \sd_{\VV_0} \big((x,t),(y,s)\big),\quad  \Sigma_2(v_1) = (1 -v_1) \sqrt{s}, \quad 
  \Sigma_3(v_2)= (1-v_2) \sqrt{1-s}.  
$$  
Recall $\zeta (x,t,y,s; v) = v_1 \sqrt{\tfrac{st + \la x,y \ra}2}+ v_2 \sqrt{1-t}\sqrt{1-s}$ and $|\zeta (x,t,y,s; v)|\le 1$. 
For simplicity, we write $\Psi = \sqrt{\frac{1+\la \xi,\eta \ra}{2}}$. Since $\la \xi,\eta \ra = \cos \sd_\SS(\xi,\eta)$
in terms of the geodesic distance on $\sph$, it follows that $\Psi = \cos \frac{\sd_{\SS}(\xi,\eta)}{2}$. Taking 
derivative of $\zeta$, 
\begin{align*} 
       \partial_t \zeta (x,t,y,s; v) =  v_1\frac{\sqrt{s}}{2 \sqrt{t}} \, \Psi  - v_2 \frac{\sqrt{1-s}}{2 \sqrt{1-t}},
\end{align*}
and writing $v_i = 1 - (1-v_i)$, it follows readily that 
$$
 |2  \varphi(t) \partial_t  \zeta (x,t,y,s; v)| \le (1 - v_1)\sqrt{s} + (1- v_2) \sqrt{1-s}+  
     \left| \sqrt{1-t} \sqrt{s}\,\Psi - \sqrt{t} \sqrt{1-s}\right|. 
$$
If $t \ge s$, we write the last term as 
$$
 \sqrt{1-t} \sqrt{s}\,\Psi - \sqrt{t} \sqrt{1-s} = - \sqrt{1-t} \sqrt{s} (1- \Psi) + \sqrt{1-t}\sqrt{s} - \sqrt{t} \sqrt{1-s},
$$
whereas if $t \le s$, we write
$$
 \sqrt{1-t} \sqrt{s}\,\Psi - \sqrt{t} \sqrt{1-s} = \left(\sqrt{1-t}\sqrt{s} - \sqrt{t}\sqrt{1-s}\right)\Psi- \sqrt{t} \sqrt{1-s} (1- \Psi).
$$
From these identities, we deduce the estimate
$$
\left|\sqrt{1-t} \sqrt{s}\,\Psi - \sqrt{t} \sqrt{1-s}\right|  \le \left| \sqrt{1-t}\sqrt{s} - \sqrt{t}\sqrt{1-s}\right| + (t s)^{\f14} (1-\Psi).
$$
Since $\Psi =  \cos \frac{\sd_{\SS}(\xi,\eta)}{2}$, we obtain an upper bound $1 - \Psi \le \f12 \sd_{\SS}(\xi,\eta)$, 
which implies by \eqref{eq:d2=d2+d2} that $(t s)^{\f14} \sd_{\SS}(\xi,\eta) \le c\,\sd_{\VV_0^{d+1}} \big( (x,t), (y,s)\big)$. 
Furthermore, we also have 
\begin{align*}
 | \sqrt{1-t} \sqrt{s}  - \sqrt{t} \sqrt{1-s} | \, & \le |\sqrt{1-t}-\sqrt{1-s}| \sqrt{s} + |\sqrt{s}-\sqrt{t}| \sqrt{1-s} \\
    &  \le 2 \sd_{\VV_0}\big( (x,t),(y,s)\big),
\end{align*}
using \eqref{eq:|s-t|}. Putting these inequalities together, we have proved \eqref{eq:Dt-zeta}. Using this estimate, 
we can then use \eqref{eq:DLn(t,1)} in \eqref{eq:DLn-intV0} to obtain for a positive number $\b$, 
\begin{align*}
 & \left |  \varphi(t) \partial_t \sL_n \big(\sw_{-1,\g}; (x,t), (y,s)\big)  \right |
   \le \,   c  \int_{[-1,1]^2} \frac{ n^{2\l+3}} 
      {\left(1+n  \sqrt{1- \zeta(x,t,y,s;v)}\right)^\beta} \\
  &  \qquad\qquad\qquad  \times  \big( \Sigma_1 + \Sigma_2(v_1) + \Sigma_3(v_2) \big)
           (1-v_1^2)^{\f{d-2}2-1}(1-v_2^2)^{\g-\f12} \d v. \notag
\end{align*}
The integral with $\Sigma_1$ term is bounded by applying Lemma \ref{lem:kernelV0} and choosing $\beta$ 
appropriately. The integral with $\Sigma_2(v_1)$ term is bounded by applying the same lemma but with 
$(1-v_1)^{\frac{d-2}{2} -1}$ replaced by $(1-v_1)^{\f{d}{2}-1}$ and then using
$n^{-1} \sqrt{s} \le (\sqrt{t} + n^{-1}) (\sqrt{s} + n^{-1})$, whereas the integral with $\Sigma_3(v_2)$ can
be handled similarly. In fact, the estimate of the last two integrals are already appeared in the proof
of \cite[Theorem 4.13]{X21}. This completes the proof.
\end{proof}

We are now ready to prove the Bernstein inequality for $\partial_t$ on the conic surface.  

\begin{thm} \label{thm:Bernstein}
Let $\sw$ be a doubling weight on $\VV_0^{d+1}$. Let $\ell$ be a positive integer and $1 \le p < \infty$.  Then
\begin{equation}\label{eq:BernV0-1}
  \|\partial_t f\|_{p,\sw} \le c n^{2\ell}  \|f\|_{p,\sw}, \quad \forall f\in \Pi_n(\VV_0^{d+1}), 
\end{equation}
and, with $\varphi(t) = \sqrt{t(1-t)}$,
\begin{equation}\label{eq:BernV0-2}
 \left \| \varphi^\ell \partial_t^\ell  f \right \|_{p,\sw} \le c n^\ell \|f\|_{p,\sw},  \quad \forall f\in \Pi_n(\VV_0^{d+1}). 
\end{equation}
\end{thm}

\begin{proof}
By the definition of $ \sL_n\left(\sw_{-1,\g}; \cdot,\cdot\right)$, every $f\in \Pi_n(\VV_0^{d+1})$ satisfies  
\begin{equation}\label{eq:reprodLn}
  f(x,t) = \int_{\VV_0^{d+1} } f(y,s) \sL_n\left(\sw_{-1,\g}; (x,t), (y,s)\right) \sw_{-1,\g}(s) \d \s(y,s).
\end{equation}
Taking $\partial_t$ derivative and using the maximal function $f_{\b,n}^\ast$ in \eqref{eq:fbn*}, we obtain
\begin{align*}
  \left |\partial_t f(x,t) \right| \le c f_{\b,n}^\ast(x,t)   \int_{\VV_0^{d+1}}
        \frac{ \left| \partial_t \sL_n\left(\sw_{-1,\g}; (x,t), (y,s)\right) \right| }
           {\big(1+n \sd_{\VV_0} ((x,t),(y,s)) \big)^\b} \sw_{-1,\g}(s) \d \sm (y,s).
\end{align*}
Using the first estimate in Lemma \ref{lem:D2kernelV0} and choosing $\b= 2 \a(\sw)/p$ and 
$\k > \b + 2 d + \frac12 (\g+\f{d-2}{2})$, we see that the integral in the right-hand side is bounded by 
$$
   \int_{\VV_0^{d+1} } |f(y,s)|  \frac{n^{d+1}  \sw_{-1,\g}(s) \d \sm(y,s)}
   {\varphi(t) \sqrt{ \sw_{\g,d} (n; t) }\sqrt{ \sw_{\g,d} (n; s) }\big(1 + n \sd_{\VV_0}( (x,t), (y,s)) \big)^{\k - \b}}
   \le c \frac{n}{\varphi(t)}
$$
where the last step follows from \eqref{eq:intLn1} with $p =1$ and $\rho = 0$. Consequently, we obtain
\begin{align}\label{eq:partial_t}
 \left| \partial_t f(x,t) \right|  \le c \frac{n}{\varphi(t)} f_{\b,n}^\ast(x,t). 
\end{align}
In particular, it follows by Proposition \ref{prop:Remz} and $\varphi(t) \ge \f1n$ for 
$\delta n^{-2} \le t \le 1- \delta n^{-2}$ that
\begin{align*}
  \|\partial_t f\|_{p,\sw} \le \|\partial_t f \chi_{n,\delta}(t)\|_{p,\sw} \le c n^2 \|f_{\b,n}^\ast \|_{p,\sw} \le  c n \|f\|_{p,\sw},
\end{align*}
where the last step follows from \eqref{eq:fbn*bound}, which proves the inequality \eqref{eq:BernV0-1} 
for $\ell =1$. Iterating the inequality for $\ell =1$ proves \eqref{eq:BernV0-1} for $\ell > 1$. In the case of 
$p=\infty$, we use $\|f\|_\infty$ instead of $f_{\b,n}^*$ in \eqref{eq:partial_t}, the resulted estimate 
immediately implies \eqref{eq:BernV0-2} and, with the help of \eqref{eq:Remz2}, also \eqref{eq:BernV0-1}. 

Since $\left| \varphi(t) \partial_t f(x,t) \right| \le c n  f_{\b,n}^\ast(x,t)$ by \eqref{eq:partial_t}, we obtain immediately, by  \eqref{eq:fbn*bound}, that $\left \|\varphi \partial_t f \right \|_{p,\sw} \le  c n \|f\|_{p,\sw}$, which proves
 \eqref{eq:BernV0-2} for $\ell =1$. Moreover, since 
$$
    \left( \varphi(t) \partial_t \right)^2 = (1-2t) \partial_t + \varphi(t)^2 \partial^2_t, 
$$
we obtain immediately that 
$$
 \left \|\varphi^2 \partial_t^2 f \right \|_{p,\sw} \le   \|\partial_t f\|_{p,\sw} + \|\left( \varphi(t) \partial_t \right)^2 f\|_{p,\sw}
  \le c n^2  \|f\|_{p,\sw},
$$
where we have applyed the inequality \eqref{eq:BernV0-1} with $\ell =1$ and the inequality \eqref{eq:BernV0-2} 
with $\ell =1$ twice in the second step. This proves the inequality \eqref{eq:BernV0-2} for $\ell =2$. Since 
$\varphi^2 \partial_t^2 f$ is a polynomial in $\Pi_n(\VV_0^{d+1}, \sw)$, we can iterate the inequality
\eqref{eq:BernV0-2} with $\ell = 2$ to establish the inequality for all even $\ell$ and, similarly, establish
the inequality for odd $\ell$ after applying the inequality with $\ell =1$ once. This completes the proof.
\end{proof}

\subsection{Bernstein inequality for $D_{i,j}$ on the conic surface}\label{sec:BernsteinV0B}
We start with the estimate of the derivative of the localized kernel. 
\begin{lem} \label{lem:DkernelV0}
Let $d\ge 2$ and $\g \ge -\f12$. Let $D_{i,j}$ be the operator acting on $x$ variable. 
Then for any $\k > 0$ and $(x,t), (y,s) \in \VV_0^{d+1}$, 
\begin{equation*}
\left | D_{i,j} \sL_n (\sw_{-1,\g}; (x,t), (y,s))\right|
\le \frac{c_\k n^{d+1} \sqrt{t}}{\sqrt{ \sw_{\g,d} (n; t) }\sqrt{ \sw_{\g,d} (n; s) }\big(1 + n \sd_{\VV_0}( (x,t), (y,s)) \big)^{\k}}
\end{equation*}
for $1 \le i< j \le d$ and, furthermore,
\begin{equation*}
\left | D_{i,j}^2 \sL_n (\sw_{-1,\g}; (x,t), (y,s))\right|
\le \frac{c_\k n^{d+2}  \big (t+ \sqrt{t}\sqrt{s}(1+ \la \xi,\eta\ra)^{-\f12}\big)}{\sqrt{ \sw_{\g,d} (n; t) }\sqrt{ \sw_{\g,d} (n; s) }
\big(1 + n \sd_{\VV_0}( (x,t), (y,s)) \big)^{\k}}. 
\end{equation*}
\end{lem}

\begin{proof} 
The inequality \eqref{eq:DLn-intV0} holds with $D_t$ replaced by $D_{i,j}$. First, we need an estimate 
of $|D_{i,j} \zeta (x,t,y,s; v)|$. If $\xi, \eta \in \sph$, then $2 \sqrt{\frac{1+\la \xi,\eta \ra}{2}} = \|\xi+\eta\|$. Hence, 
for $x = t\xi$ and $y = s \eta$, 
$$
  D_{i,j} \zeta (x,t,y,s; v) = v_1  \frac{x_i y_j - x_j y_i} {4\sqrt{\tfrac{t s + \la x, y\ra}2}}
     = v_1\sqrt{t}\sqrt{s}\,\frac{\xi_i\eta_j - \xi_j \eta_i} {2 \|\xi + \eta\|}.
$$
Now, using the identity
\begin{align*}
  2 (\xi_i \eta_j - \xi_j \eta_i) = (\xi_i - \eta_i ) (\xi_j + \eta_j ) - (\xi_j - \eta_j ) (\xi_i + \eta_i ) 
\end{align*} 
and $\|\xi - \eta\| \le \sd_{\SS}(\xi,\eta)$ for $\xi,\eta \in \sph$, it follows immediately that 
$$
   | \xi_j \eta_j - \xi_j \eta_i | \le \|\xi- \eta\|\, \| \xi + \eta \| \le  \|\xi+ \eta\|\,\sd_{\SS}(\xi,\eta).
$$
Putting these together, we obtain the inequality  
\begin{align*}
   | D_{i,j}  \zeta (x,t,y,s; v)| \le |v_1| \sqrt{t}\sqrt{s} \, \sd_{\SS} (\xi,\eta) \le \sqrt{t}\sqrt{s} \, \sd_{\SS} (\xi,\eta). 
\end{align*}
If $s \le t$, then $\sqrt{s} \, \sd_{\SS} (\xi,\eta) \le (t s)^\f14  \sd_{\SS} (\xi,\eta) \le c\, \sd_{\VV_0}((x,t), (y,s))$ by \eqref{eq:d2=d2+d2}. If $s \ge t$, then $\sqrt{s} \le (t s)^{\f14} + \sd_{\VV_0}((x,t), (y,s))$ by \eqref{eq:|s-t|}, 
which implies, by \eqref{eq:d2=d2+d2},
$$
\sqrt{s} \,\sd_{\SS} (\xi,\eta)  \le (t s)^{\f14}  \sd_{\SS} (\xi,\eta) + \pi \sd_{\VV_0}((x,t), (y,s)) \le c \, \sd_{\VV_0}((x,t), (y,s)).
$$
Consequently, we have proved the inequality 
\begin{align}\label{eq:Dijz-bd}
   | D_{i,j}  \zeta (x,t,y,s; v)| \le c \sqrt{t} \sd_{\VV_0}((x,t), (y,s)). 
\end{align}
Using this inequality and $|\zeta (x,t,y,s;v)| \le 1$, we then use \eqref{eq:DLn(t,1)} in \eqref{eq:DLn-intV0} to obtain 
\begin{align*}
  \left | D_{i,j} \sL_n \big(\sw_{-1,\g}; (x,t), (y,s)\big)  \right |
     \le \, & c  \int_{[-1,1]^2} \frac{ n^{2\l+3}  \sqrt{t} \sd_{\VV_0}\big((x,t),(y,s)\big)} 
      {\left(1+n  \sqrt{1- \zeta(x,t,y,s;v)}\right)^\k} \\
  &    \times  (1-v_1^2)^{\f{d-2}2-1}(1-v_2^2)^{\g-\f12} \d v \notag
\end{align*}
from which the desired estimate follows from Lemma \ref{lem:kernelV0}. 

Now, taking the derivative $D_{i,j}$ one more time, we obtain 
\begin{align*}
   &D_{i,j}^2 \sL_n (\sw_{-1,\g}; (x,t), (y,s) ) \\
   = \, &   c_{\g}  \int_{[-1,1]^2}  \partial^2 L_n ^{(\l-\f12,-\f12)}\big(2 \zeta (v)^2-1 \big)
   \left(4  \zeta (v) D_{i,j} \zeta (v) \right)^2 (1-v_1^2)^{\f{d-2}2-1}(1-v_2^2)^{\g-\f12} \d v \\
   + \,&   c_{\g}  \int_{[-1,1]^2}  \partial L_n ^{(\l-\f12,-\f12)}\big(2 \zeta (v)^2-1 \big)
       4 D_{i,j} \left[ \zeta (v) D_{i,j} \zeta (v) \right](1-v_1^2)^{\f{d-2}2-1}(1-v_2^2)^{\g-\f12} \d v.
\end{align*}
Applying \eqref{eq:DLn(t,1)} with $m=2$ and \eqref{eq:Dijz-bd}, it is easy to see that the first integral
in the right-hand side is bounded by 
\begin{align*}
 c  \int_{[-1,1]^2} \frac{ n^{2\l+5}  t \sd_{\VV_0}\big((x,t),(y,s)\big)^2} 
      {\left(1+n  \sqrt{1- \zeta(x,t,y,s;v)}\right)^\k} (1-v_1^2)^{\f{d-2}2-1}(1-v_2^2)^{\g-\f12} \d v,
\end{align*}
which can be estimated as in the case of $D_{i,j} L_n$ to give the desired upper bound. Moreover,
the second integral is a sum of two integrals according to 
$$
   D_{i,j} \left[ \zeta (v) D_{i,j} \zeta (v) \right] =  \left[D_{i,j} \zeta (v) \right]^2 + \zeta(v) D_{i,j}^2\zeta(v);
$$
while the first one, using \eqref{eq:Dijz-bd}, is dominated by the first integral, the second one 
can be bounded using the identity 
$$
   D_{i,j}^2 \zeta(x,t,y,s;v) = - v_1 \sqrt{t}\sqrt{s}  \left( \frac{\xi_i \eta_i + \xi_j \eta_j}{2 \|\xi+\eta\|} 
    + \frac{ (\xi_i \eta_j- \xi_j \eta_i)^2}{2 \|\xi+\eta\|^3}\right),
$$
which implies immediately that 
$$
   |D_{i,j}^2 \zeta(x,t,y,s;v)| \le  \sqrt{t}\sqrt{s} \frac{1+\sd_{\SS}(\xi,\eta)^2}{2 \|\xi+\eta\|} \le 
     c \frac{\sqrt{t}\sqrt{s}}{\sqrt{1+\la \xi,\eta\ra}}
$$
and leads to the upper bound 
$$
   c  \frac{\sqrt{t}\sqrt{s}}{\sqrt{1+\la \xi,\eta\ra}} \int_{[-1,1]^2} \frac{ n^{2\l+3}} 
      {\left(1+n  \sqrt{1- \zeta(x,t,y,s;v)}\right)^\k} (1-v_1^2)^{\f{d-2}2-1}(1-v_2^2)^{\g-\f12} \d v
$$
which is again seen to be bounded by the desired upper bound. This completes the proof. 
\end{proof}

\begin{thm} \label{thm:Bernstein2}
Let $\sw$ be a doubling weight on $\VV_0^{d+1}$. If $f\in \Pi_n(\VV_0^{d+1})$ and $1\le p < \infty$,
then for $1 \le i,j \le d$, 
\begin{equation}\label{eq:BernV0-3}
    \left \|   \frac{1}{\sqrt{t}} D_{i,j} f \right \|_{p,\sw} \le c n \|f\|_{p,\sw} 
         \quad \hbox{and}\quad   \left \| \frac{1}{t} D_{i,j}^2 f \right \|_{p,\sw} \le c n^2 \|f\|_{p,\sw}.
\end{equation}
\end{thm}

\begin{proof}
For $f\in \Pi_n(\VV_0^{d+1})$, we start with the identity \eqref{eq:reprodLn}. Taking $D_{i,j}$ derivative 
and using the maximal function $f_{\b,n}^\ast$ in \eqref{eq:fbn*}, we obtain 
\begin{equation} \label{eq:Dijf-bd}
 \left| D_{i,j}^\ell f(x,t) \right| \le  f_{\b,n}^\ast(x,t) \int_{\VV_0^{d+1} } 
 \frac{ \left| D_{i,j}^\ell \sL^n(\sw_{-1,\g}; (x,t),(y,s)) \right| }{\big(1+n\sd_{\VV_0}((x,t),(y,s))\big)^\b} 
 \sw_{-1,\g}(s) \d \sm (y,s).
\end{equation}
Using the estimate in Lemma \ref{lem:DkernelV0} and choosing $\b= 2 \a(\sw)/p$ and $\k > \b + 2 d + 
\frac12 (\g+\f{d-2}{2})$, it follows that, for $\ell =1$ the integral in the right-hand side is bounded by 
$$
c  \int_{\VV_0^{d+1} } \frac{n^{d+1} \sqrt{t}  \sw_{-1,\g}(y,s) \d \sm (y,s)}
   {\sqrt{ \sw_{\g,d} (n; t) }\sqrt{ \sw_{\g,d} (n; s) }\big(1 + n \sd_{\VV_0}( (x,t), (y,s)) \big)^{\k-\b}} \le c n \sqrt{t},
$$
where  the last step follows from \eqref{eq:intLn1} with $p =1$ and $\rho = 0$. Consequently, we obtain
\begin{align*}
 \left| \frac{1}{\sqrt{t}}D_{i,j} f(x,t) \right| \le c n  f_{\b,n}^\ast(x,t). 
 \end{align*}
For $1 \le p < \infty$, we can then apply \eqref{eq:fbn*bound} for $\VV_0^{d+1}$ to conclude that
$$
\left\| \frac{1}{\sqrt{t}}D_{i,j} f\right\|_{p,\sw} \le c n \| f_{\b,n}^\ast\|_{p,\sw} \le c n \| f \|_{p,\sw}.
$$
This proves \eqref{eq:BernV0-3} for $D_{i,j}$. For $D_{i,j}^2$, we first observe that, using
$$
   D_{i,j}^2  = - x_i \partial_i - x_j \partial_j + x_i D_{i,j} \partial_j - x_j D_{i,j} \partial_i,
$$
the function $t^{-1} D_{i,j}^2 f(t\xi, t)$ is a polynomial in the variables $t$ and $\xi$. Hence, by Proposition 
\ref{prop:Remz} and Remark \ref{rem:remark1}, there is a $\delta> 0$ such that, with 
$\chi_{n,\d}(t) = \chi_{[\f{\delta}{n}, 1- \f{\delta}{n}]}(t)$, 
$$
 \left \| \frac{1}{t} D_{i,j}^2 f \right \|_{p,\sw} \le c_\delta  \left \| \frac{1}{t} D_{i,j}^2 f \cdot \chi_{n,\delta} \right \|_{p,\sw}.
$$
Using $\sqrt{s} \le \sqrt{t} + \sd_{\VV_0}((x,t),(y,s))$ and the estimate in Lemma \ref{lem:DkernelV0}, the
integral in the right-hand side of \eqref{eq:Dijf-bd} with $\ell =2$ is bounded by 
\begin{align*}
  &  c_\k  t \int_{\VV_0^{d+1} }  \frac{n^{d+2}(1+\la \xi,\eta\ra)^{-\f12}
      \sw_{-1,\g}(y,s) \d \sm (y,s)}
   {\sqrt{ \sw_{\g,d} (n; t) }\sqrt{ \sw_{\g,d} (n; s) }\big(1 + n \sd_{\VV_0}( (x,t), (y,s)) \big)^{\k-\b}} \\
   & + c_\k t \int_{\VV_0^{d+1} } \frac{n^{d+1} (1+\la \xi,\eta\ra)^{-\f12} \sw_{-1,\g}(y,s) \d \sm (y,s)}
   {\sqrt{t} \sqrt{ \sw_{\g,d} (n; t) }\sqrt{ \sw_{\g,d} (n; s) }\big(1 + n \sd_{\VV_0}( (x,t), (y,s)) \big)^{\k-\b-1}}, 
\end{align*}
which, for $t \in [\f{\delta}{n}, 1- \f{\delta}{n}]$ so that $1/\sqrt{t} \le n$, is further bounded by 
$$
 c_\k t \int_{\VV_0^{d+1} } |f(y,s)|  \frac{n^{d+2} (1+\la \xi,\eta\ra)^{-\f12} 
   \sw_{-1,\g}(y,s) \d \sm (y,s)}
   {\sqrt{t} \sqrt{ \sw_{\g,d} (n; t) }\sqrt{ \sw_{\g,d} (n; s) }\big(1 + n \sd_{\VV_0}( (x,t), (y,s)) \big)^{\k-\b-1}}
  \le c_\k t\, n^2,  
$$
by applying \eqref{eq:intLn1} with $p =1$ and $\rho = \f12$. Thus, we obtain  
$$
 \left| \frac{1}{t}D_{i,j}^2 f(x,t) \right| \chi_{n,\delta}(t) \le c n^2  f_{\b,n}^\ast(x,t),
$$
from which the proof follows readily from the boundedness of $f_{\b,n}^*$. Using 
\eqref{eq:Remz2} and $\|f\|_\infty$ instead of $f_{\b,n^*}$, the above proof also applies to the 
case $p = \infty$.
\end{proof}

The inequalities in \eqref{eq:BernV0-3} are the cases of $\ell = 1$ and $\ell =2$ of \eqref{eq:BI-V0-3}. We
now show that the inequality \eqref{eq:BI-V0-3} does not hold in general for $\ell > 2$, we consider for
example the polynomial $f_0(x,t) = x_1$. It is easy to see that 
$$
  D_{1,j}^{2\ell} f_0(x,t) = (-1)^\ell x_1 \quad \hbox{and}\quad    D_{1,j}^{2\ell -1} f_0(x,t) = (-1)^\ell x_j,
$$
so that $(\sqrt{t})^{-\ell} D_{i,j}^\ell f_0(t\xi,t)$ has a singularity $t^{-\ell/2 +1}$ for $\ell > 2$. Hence, taking
$\sw = \sw_{-1,\g}$ on $\VV_0^3$, for example, it is easy to see that $\|t^{-\ell/2} D_{i,j}^\ell f_0\|_{2,\sw_{-1,\g}}$ 
is infinite if $\ell \ge  3$.  Thus, \eqref{eq:BI-V0-3} does not hold for $\ell \ge 3$ in general. 

Finally, for $\ell = 1,2, \ldots$ and any doubling weight $\sw$, the inequality \eqref{eq:BI-V0-4} holds. 
Indeed, the case $\ell = 1$ of \eqref{eq:BI-V0-4} follows readily from, and in fact weaker than, 
\eqref{eq:BernV0-3}. The case $\ell > 1$ of \eqref{eq:BI-V0-4} follows readily from iteration. 

\subsection{Self-adjoint of $-\Delta_{-1,\g}$} \label{sect:proof_corollary}
In this subsection we give a proof of Proposition \ref{prop:self-adJacobi} and Corollary \ref{cor:B-V0_p=2}. 

\medskip
\begin{proof}[Proof of Proposition \ref{prop:self-adJacobi}] 
First we assume $\g > 0$. The differential operator can be written as 
$$
    \Delta_{0,\g} = \frac{1}{t^{d-2}(1-t)^\g} \frac{\partial}{\partial t} \left ( t^{d-1} (1-t)^{\g+1}  \frac{\partial}{\partial t}  \right)
    + t^{-1} \Delta_0^{(\xi)}. 
$$
Setting $x = t \xi$, the integral over the conic surface can be written as 
$$
\int_{\VV_0^{d+1}}f(x,t)  \d\sm(x,t) = \int_0^1 t^{d-1} \int_{\sph} f(t \xi, t) \d \s(\xi) \d t. 
$$
Hence, we can write 
\begin{align*}
 \int_{\VV_0^{d+1}} & \Delta_{0,\g} f(x,t) \cdot g(x,t) \sw_{-1,\g}(t)  \d\sm(x,t) \\
 = & \int_0^1 \int_{\sph}   \frac{\partial}{\partial t} \left ( t^{d-1} (1-t)^{\g+1} \frac{\partial f}{\partial t}\right)
      \cdot g(t \xi, t) \d \s(\xi)\d t \\
    & + \sum_{1 \le i < j \le d} \int_0^1  t^{d-2}  \int_{\sph}  D_{i,j}^2 f(t \xi,t) \cdot g(t \xi, t) \d \s(\xi)  \sw_{-1,\g}(t) \d t.
\end{align*}
We then integrate by parts in $t$ variable, using the assumption that $\g > 0$ and $d \ge 2$, in
the first integral in the right-hand side, and applying (cf. \cite[Proposition 1.8.4]{DaiX})
\begin{equation}\label{eq:Dij-selfadj}
   \int_{\sph} f(\xi)D_{i,j} g(\xi) \d \s(\xi) = - \int_{\sph} D_{i,j}f(\xi) g(\xi) \d \s(\xi), \quad 1 \le i,j \le d,
\end{equation}
on the second integral in the right-hand side to establish the desired formula. Analytic continuation shows that the 
identity holds for $\g > -1$. 
\end{proof}

\begin{proof}[Proof of Corollary \ref{cor:B-V0_p=2}]
The inequality follows from setting $g =f$ in the Proposition \ref{prop:self-adJacobi} and the Cauchy inequality.
If $f \in \Pi_n$, then the inequality \eqref{eq:BernsteinLB} shows that the right-hand side of the stated 
inequality is bounded by $c n^2 \left \| f \right \|_{\sw, 2}^2$, so that each term of the left-hand sider is bounded
by the same quantity and, hence, we obtain \eqref{eq:BI-V0-2} and \eqref{eq:BI-V0-3} for $\ell =1$ and $p=2$.
\end{proof}

\section{Bernstein inequalities on the cone}\label{sec:OP-cone}
\setcounter{equation}{0}

We work with the solid cone $\VV^{d+1}$ in this section. In the first subsection, we recall the orthogonal 
polynomials and their kernels on the cone. Instead of estimating the derivatives of the kernels as in the 
previous section, we rely on an intimate connection between analysis on the cone and on the conic surface,
which is explained in the second subsection. The proof of the Bernstein inequalities on the cone, based 
on the connection, is given in the third subsection. Finally, the proof of Theorem \ref{thm:self-adjointV} 
and its corollary is in the fourth subsection. 

\subsection{Orthogonal polynomials and localized kernels}\label{sec:OPcone}
Let $\Pi_n^d$ denote the space of polynomials of degree at most $n$ in $d$ variables. For $\mu > -\f12$ 
and $\g > -1$, we define the weight function $W_{\mu,\g}$ by
$$
   W_{\mu,\g}(x,t) = (t^2-\|x\|^2)^{\mu-\f12} (1-t)^\g, \qquad \|x\| \le t, \quad 0 \le t \le 1. 
$$
Orthogonal polynomials with respect to $W_{\mu,\g}$ on $\VV^{d+1}$ are studied in \cite{X20}, which are 
orthogonal with respect to the inner product 
$$
\la f, g\ra_{W_{\mu,\g}} =\bs_{\mu,\g} \int_{\VV^{d+1}} f(x,t) g(x,t) W_{\mu,\g} \d x \d t.
$$ 
Let $\CV_n(\VV^{d+1},W_{\mu,\g})$ be the space of orthogonal polynomials of degree $n$. Then 
$$
   \dim \CV_n(\VV^{d+1},W_{\mu,\g})  = \binom{n+d+1}{n}, \quad n=0, 1,2,\ldots.
$$
An orthogonal basis of $\CV_n(\VV^{d+1}, W_{\mu,\g})$ can be given in terms of the Jacobi polynomials 
and the orthogonal polynomials on the unit ball. For $m =0,1,2,\ldots$, let $\{P_\kb^m(W_\mu): |\kb| = m, 
\kb \in \NN_0^d\}$ be an orthonormal basis of $\CV_n(\BB^d, W_\mu)$ on the unit ball. Let
\begin{equation} \label{eq:coneJ}
  \Jb_{m,\kb}^n(x,t):= P_{n-m}^{(2m+2\mu+d-1, \g)}(1- 2t) t^m P_\kb^m\left(W_\mu; \frac{x}{t}\right), \quad 
\end{equation}
Then $\{\Jb_{m,\kb}^n(x,t): |\kb| = m, \, 0 \le m\le n\}$ is an orthogonal basis of $\CV_n(\VV^{d+1},W_{\mu,\g})$.
The following theorem is established in \cite[Theorem 3.2]{X20}. 

\begin{thm}\label{thm:Delta0V}
Let $\mu > -\tfrac12$, $\g > -1$ and $n \in \NN_0$. Let $\fD_{\mu,\g}$ be the second order differential 
operator defined in \eqref{eq:V-DE}. Then the polynomials in $\CV_n(\VV^{d+1},W_{\mu,\g})$ are 
eigenfunctions of $\fD_{\mu,\g}$; more precisely,  
\begin{equation}\label{eq:cone-eigen}
   \fD_{\mu,\g} u =  -n (n+2\mu+\g+d) u, \qquad \forall u \in \CV_n(\VV^{d+1},W_{\mu,\g}).
\end{equation}
\end{thm} 
 
The reproducing kernel of the space $\CV_n(\VV^{d+1}, W_{\mu,\g})$, denoted by $\Pb_n(W_{\g,\mu};\cdot,\cdot)$, 
can be written in terms of the above basis, 
$$
\Pb_n\big(W_{\mu,\g}; (x,t),(y,s) \big) = \sum_{m=0}^n \sum_{|\kb|=n} 
    \frac{  \Jb_{m, \kb}^n(x,t)  \Jb_{m, \kb}^n(y,s)}{\la  \Jb_{m, \kb}^n,  \Jb_{m, \kb}^n \ra_{W_{\mu,\g}}}.
$$
It is the kernel of the projection $\proj_n(W_{\mu,\g}): L^2(\VV^{d+1},W_{\mu,\g}) \to 
\CV_n(\VV^{d+1}, W_{\mu,\g})$, 
$$
\proj_n(W_{\mu,\g};f) = \int_{\VV^{d+1}} f(y,s) \Pb_n\big(W_{\mu,\g}; \,\cdot, (y,s) \big)  W_{\mu,\g}(s) \d y \d s.
$$
This kernel enjoys an addition formula that can be written as a triple integral over $[-1,1]^3$ for $\mu > 0$; see 
\cite[Theorem 4.3]{X20}. For our purpose, we only need the closed-form formula for the degenerate case 
$\mu = 0$; see \cite[Theorem 5.8]{X21}.

\begin{thm} \label{thm:PnCone2}
Let $d \ge 2$, $\mu =0$ and $\g \ge -\f12$. Then, for $n =0,1,2,\ldots$, 
\begin{align}\label{eq:PbCone3}
   \Pb_n \big(W_{0,\g};  (x,t), (y,s)\big) =  c_{\g,d}
  &  \int_{[-1,1]^2}   \left[ Z_{n}^{(\g+d-\f12,-\f12)} (\xi (x, t, y, s; 1, \vb))  \right. \\
  & \qquad\quad \left. +Z_n^{(\g+d-\f12,-\f12)}(\xi (x, t, y, s; -1, \vb))\right] \notag \\
 & \times (1-v_1^2)^{\f{d-1}2 - 1}(1-v_2^2)^{\g-\f12}  \d \vb, \notag
\end{align}
where $c_{\g,d}$ is a constant, so that $\Pb_0 =1$ and $\xi (x,t, y,s; u, \vb) \in [-1,1]$ is defined by 
\begin{align} \label{eq:xi}
\xi (x,t, y,s; u, \vb) = &\, v_1 \sqrt{\tfrac12 \left(ts+\la x,y \ra + \sqrt{t^2-\|x\|^2} \sqrt{s^2-\|y\|^2} \, u \right)}\\
      & + v_2 \sqrt{1-t}\sqrt{1-s}. \notag
\end{align} 
\end{thm}

Let $\wh a$ be an admissible cut-off function.  For $(x,t)$, $(y,s) \in \VV^{d+1}$, the localized kernel 
$\Lb_n(W_{0,\g}; \cdot,\cdot)$ is defined by  
$$
   \Lb_n\left(W_{0,\g}; (x,t),(y,s)\right) = \sum_{j=0}^\infty \wh a\left( \frac{j}{n} \right)
       \Pb_j\left(W_{0,\g}; (x,t), (y,s)\right). 
$$
Like the kernel $\sL_n$ on the conic surface, the kernel $\Lb_n(W_{0,\g})$ is also highly localized via
the distance function $\sd_{\VV}(\cdot, \cdot)$ on the cone. There is, however, no need to carry out
the estimate for these kernels as in the case of the conic surface. Indeed, there is a direct connection between
the kernels for $W_{0,\g}$ on $\VV^{d+1}$ and those for $w_{-1,\g}$ on the conic surface $\VV_0^{d+2}$ 
in $\RR^{d+2}$, as will be explained in the next subsection.

\subsection{The analysis on $\VV^{d+1}$ and on $\VV_0^{d+2}$} 
Let $(x,t)$ and $(y,s)$ be the elements of the solid cone $\VV^{d+1}$. We shall adopt the notation 
\begin{equation} \label{eq:XY}
 X:=\big(x,\sqrt{t^2-\|x\|^2}\big), \quad Y:=\big(y,\sqrt{s^2-\|y\|^2}\big), \quad Y_*:=\big(y, - \sqrt{s^2-\|y\|^2}\big).
\end{equation}
It is evident that $\|X\| = t$, $\|Y\| = \|Y_*\| =s$, so that $(X,t)$, $(Y,s)$ and $(Y_*,s)$ are elements 
of the conic surface $\VV_0^{d+2}$ in $\RR^{d+2}$. 

\begin{prop}
Let $\sL_n(w_{-1,\g}; \cdot,\cdot)$ be the kernel defined in \eqref{eq:Ln-intV0} for $\VV_0^{d+2}$. Then
\begin{equation} \label{eq:LnW-Lnw}
   \Lb_n \left(W_{0,\g}; (x,t), (y,s)\right) = 
       \sL_n \big(\sw_{-1,\g}; (X,t), (Y,s)\big) + \sL_n \big(\sw_{-1,\g}; (X,t), (Y_*,s)\big).
\end{equation}
\end{prop}

\begin{proof}
By \eqref{eq:PbCone3}, we can write 
$$
 \Lb_n\left(W_{0,\g}; (x,t), (y,s)\right) = \Lb_n^+ \left(W_{0,\g}; (x,t), (y,s)\right)+ \Lb_n^-\left(W_{0,\g}; (x,t), (y,s)\right),
$$
where, with $\l = \g+d$, 
\begin{align*}
\Lb_n^\pm\left(W_{0,\g}; (x,t), (y,s)\right) = c_{0,\g,d} \int_{[-1,1]^2}&L_n^{(\l-\f12, -\f12)}\left(2 \xi(x,t,y,s; \pm 1,v)^2-1\right)\\
     &  \qquad  \times (1-v_1^2)^{\f{d-1}{2}-1}(1-v_2^2)^{\g-\f12} \d v.
\end{align*}
Let $X = t X'$ and $Y = s Y'$. Then $X' = (x',\sqrt{1-\|x'\|^2}) \in \SS^d$ and $Y' = (y',\sqrt{1-\|y'\|^2}) \in \SS^d$. 
Then, using the explicit formula of $\xi$ in \eqref{eq:xi} and $\zeta$ in \eqref{eq:zetaV0}, it follows
\begin{equation*}
   \xi(x,t,y,s; 1,v) = v_1 \sqrt{t s} \sqrt{\tfrac12 (1+ \la X',Y'\ra)}+v_2 \sqrt{1-t}\sqrt{1-s} 
        = \zeta(X ,t, Y, s; 1,v). 
\end{equation*}
The similar identity holds for $\xi(x,t, y,s; -1, v)$ if we replace $Y$ by $Y_*$ in the right-hand side. As a consequence,
comparing with \eqref{eq:Ln-intV0}, we obtain
\begin{align*}
         \Lb_n^+ \left(W_{0,\g}; (x,t), (y,s)\right) &= \sL_n \left(\sw_{-1,\g}; (X,t), (Y,s)\right), \\
         \Lb_n^- \left(W_{0,\g}; (x,t), (y,s)\right) & = \sL_n \left(\sw_{-1,\g}; (X,t), (Y_*,s)\right),
\end{align*}
from which the stated identity follows readily. 
\end{proof}

The connection between the cone $\VV^{d+1}$ and the conic surface $\VV_0^{d+2}$ is further manifested in
the following integral relation. 

\begin{prop}\label{prop:IntV0V}
Let $f: \RR^{d+1} \mapsto \RR$ be a continuous function. Let $x_{d+1} = \sqrt{t^2-\|x\|^2}$ for $\|x\| \le t$. Then
$$
  \int_{\VV_0^{d+2}} f(y,s) w(s) \d \sm(y,s) = \int_{\VV^{d+1}} 
    \big [ f\big((x,x_{d+1}),t\big) +  f\big((x, - x_{d+1}),t\big)\big ] t w(t) \frac{\d x\d t}{\sqrt{t^2-\|x\|^2}} .
$$
\end{prop}

\begin{proof}
Let $\d \s$ be the surface measure on $\SS^d$. We write the integral on $\VV_0^{d+1}$ as
\begin{align*}
&\int_{\VV_0^{d+2}} f(y,s) w(s) \d \sm(y,s) = \int_0^1 s^d  \int_{\SS^d} f(s \xi,s) w(s) \d\s(\xi) \d s \\
& = \int_0^1 s^d \int_{\BB^d} 
   \left[ f\big(s (u,\sqrt{1-\|u\|^2}), s\big) +  f\big(s (u, - \sqrt{1-\|u\|^2}), s\big) \right] 
      \frac{\d u}{\sqrt{1-\|u\|^2}} \d x w(s) \d \s,
\end{align*}
where the second identity uses \cite[(A.5.4)]{DaiX}. Rewriting the right-hand side as the integral over
$\VV^{d+1}$, with $x = s u$, proves the stated identity.
\end{proof}

The doubling weight on $\VV^{d+1}$ is defined as usual via the distant function $\sd_{\VV}$. The latter 
can be defined in terms of the distance function on the conic surface $\VV_0^{d+2}$. Indeed, using 
$X$ and $Y$ in \eqref{eq:XY}, we define
\begin{equation} \label{eq:distV}
  \sd_{\VV} ((x,t), (y,s)) : =  \arccos \left(\sqrt{\frac{\la X,Y\ra +ts}2}  + \sqrt{1-t}\sqrt{1-s}\right),
\end{equation}
which is indeed a distance function on the solid cone $\VV^{d+1}$, as shown in \cite{X21}. In particular, 
we have the relation 
\begin{equation}\label{eq:dV=dV0}
    \sd_{\VV^{d+1}} ((x,t), (y,s)) =  \sd_{\VV_0^{d+2}} ((X,t), (Y,s)), \quad (x,t), (y,s) \in \VV^{d+1},
\end{equation}
where we write $\sd_{\VV} =  \sd_{\VV^{d+1}}$ and $ \sd_{\VV_0} = \sd_{\VV_0^{d+2}}$ to emphasis
the dimension of the domains.

\subsection{Proof of main results} 
We will need an analog of Proposition \ref{prop:Remz} that is of interest in itself.  
 
\begin{prop} \label{prop:RemzV}
Let $W$ be a doubling weight function on $\VV^{d+1}$. For $n \in \NN$, let $\chi_{n,\delta}(x,t)$ denote the 
characteristic function of the set 
$$
   \left \{(x,t): \frac \delta{n^2} \le t \le 1- \frac \delta {n^2}, \quad \sqrt{t^2- \|x\|^2} \ge \delta \frac{\sqrt{t}}{n} \right \}.
$$
Then, for $f \in \Pi_n^d$, $1 \le p < \infty$, and every $\delta > 0$, 
\begin{equation}\label{eq:RemzV}
  \int_{\VV^{d+1}} |f(x,t)|^p W(x,t) \d x \d t \le c_\delta \int_{\VV^{d+1}} |f(x,t)|^p \chi_{n,\delta}(x,t) W(x,t) \d  x \d t.
\end{equation}
\end{prop}

\begin{proof}
The proof uses the Marcinkiewicz-Zygmund inequality on the cone established in \cite[Theorem 5.6]{X21} 
that uses a maximal $\ve$-separated subset $\Xi_{\VV}$ on $\VV^{d+1}$. It is almost verbatim as the proof 
of Proposition \ref{prop:Remz}, once it is shown that if $(x,t) \in \Xi_{\VV}$ and $\ve \approx n^{-1}$, then 
\begin{equation} \label{eq:tj-epjV}
   \frac \delta{n^2} \le t \le 1- \frac \delta {n^2} \quad\hbox{and}\quad \sqrt{t^2- \|x\|^2} \ge \delta \frac{\sqrt{t}}{n}.
\end{equation}
The set $\Xi_\VV$ is constructed in \cite[Proposition 5.15]{X21} using the separated sets in the 
$t$-variable on $[0,1]$ and in the $x'$ variable on $\BB^d$, where $(x,t) = (t x', t) \in \VV^{d+1}$. The 
construction is similar to that of $\Xi_{\VV_0}$. In fact, for $\ve > 0$ and $N = \lfloor \frac{\pi}{2}\ve^{-1} \rfloor$,
we choose the same $t_j$ and $\ve_j$ as given in \eqref{eq:tj-epj} and define 
$$
       \Xi_{\VV} = \big\{(t_j x', t_j): \,   x \in \Xi_\BB(\ve_j), \, 1\le j \le N \big\},
$$
where $\Xi_\BB(\ve_j)$ is the maximal $\ve_j$-separated set of $\BB^d$. Choose again $\ve = \b n^{-1}$, 
$\b > 0$, so that $\ve_j \sim \sqrt{t_j} /n$. By the choices of $t_j$, the bounds for $t_j$ in \eqref{eq:tj-epjV} is 
immediate. Setting $x = t x'$, $x' \in \BB^d$, so that $\sqrt{t^2 - \|x\|^2} = t \sqrt{1-\|x'\|^2}$, we see that the 
second inequality in \eqref{eq:tj-epjV} is equivalent to, since $\ve_j = \frac{\pi}{2} \frac{\ve}{\sqrt{t_j}}$, 
$$
   \sqrt{1-\|x'\|^2} \ge c\, \ve_j, \qquad x' \in \Xi_{\BB}(\ve_j), \quad 1 \le j \le n.
$$
Thus, the question is reduced to the existence of a maximal $\ve$-separated set of $\BB^d$ that satisfies 
the above inequality. A construction is given below. 

We show that there is a maximal $\ve$-separated set $\Xi_\BB(\ve)$ on $\BB^d$ such that 
\begin{equation} \label{eq:Xi_B}
   \sqrt{1-\|u\|^2} \ge c \ve, \qquad u \in \Xi_\BB(\ve) 
\end{equation}
for any given $\ve> 0$. We choose $N =  \lfloor \frac{\pi}{2}\ve^{-1} \rfloor$ as above. For $1\le j \le N$ 
we define 
$$
 \t_j:= \frac{(2 j-1)\pi}{2 N},  \qquad \t_j^- :=  \t_j- \frac{\pi}{2 N}  \quad \hbox{and} \quad \t_j^+ :=  \t_j +\frac{\pi}{2 N}.
$$
Let $r_j =  \sin \frac{\t_j}2$ and define $r_j^-$ and $r_j^+$ accordingly. In particular, $r_1^- = 0$ and $r_N^+ = 1$. 
Then $\t_{j+1}^- =\t_j^+$ and $\BB^d$ can be partitioned by 
$$
   \BB^d =  \bigcup_{j=1}^N \BB_0^{(j)}, \quad \hbox{where}\quad \BB_0^{(j)}:= 
        \left\{ r \xi \in \BB^d:   r_j^- < r \le r_j^+, \, \xi \in \sph \right\}.  
$$
Let $\s_j := (2 r_j)^{-1} \pi \ve$.  Let $\Xi_\SS(\s_j)$ be the maximal $\s_j$-separated set of $\sph$,
such that $\{\SS_\xi(\s_j): \xi \in \Xi_\SS(\s_j)\}$ is a partition of the sphere, 
$\sph = \bigcup_{\eta \in \Xi_\SS(\s_j)}\SS_\eta(\s_j)$. We define 
$$
 \Xi_{\BB}(\ve) = \left\{ r_j \xi: \xi \in \Xi_\SS(\s_j), \, 1 \le j \le N \right\}. 
$$ 
By Definition \ref{defn:separated-pts}, to show that $\Xi_\BB(\ve)$ is $\ve$-separated, we need 
$\sd_\BB(r_j \xi, r_k \eta) \ge \ve$ for any two distinct points $r_j \xi$ and $r_k \eta$ in $\Xi_\BB(\ve)$.
It is known that the distance function $\sd_\BB$ of $\BB^d$ is defined by 
$$
  \sd_{\BB}(x,y) = \arccos \left( \la x,y\ra + \sqrt{1-\|x\|^2} \sqrt{1-\|y\|^2} \right), \qquad x, y \in \BB^d.
$$
Let $x = t \xi$ and $y = s \eta$, for $t,s \in [0,1]$ and $\xi,\eta \in \sph$, it is easy to verify that 
\begin{equation} \label{eq:dB}
   1 - \cos \sd_{\BB}(x,y) = 1- \cos \sd_{[0,1]}(t^2,s^2) + t s \left(1-\cos \sd_{\SS}(\xi,\eta) \right) 
\end{equation}
using $\cos \sd_{[0,1]}(t,s) = \sqrt{t}\sqrt{s} + \sqrt{1-t}\sqrt{1-s}$. 
We use \eqref{eq:dB} with $1- \cos \phi = 2 \sin^2 \f\phi 2$ and also use $\sd_{[0,1]}(t,s) = \f12 |\t - \phi|$ if 
$t = \sin^2 \f{\t}{2}$ and $s = \sin^2 \f{\phi}2$. If $j \ne k$ then, by \eqref{eq:dB}, 
$$
  \sd_\BB(r_j \xi, r_k \eta) \ge  \sd_{[0,1]}(r_j^2, r_k^2) =\frac12 |\t_j  - \t_k| \ge \frac{\pi}{2N} \ge \ve. 
$$
If $j = k$, then $\xi$ and $\eta$ belong to the same $\Xi_{\SS}(\s_j)$, so that $\sd_{\SS}(\xi,\eta) \ge \s_j$. 
Hence, using $\frac{2}{\pi} \phi \le \sin \phi \le \phi$, we deduce from \eqref{eq:dB} that 
$$
    \sd_\BB(r_j \xi,r_j\eta) \ge\frac{2}{\pi} r_j \sd_\SS(\xi,\eta) \ge \frac{2}{\pi} r_j \s_j = \ve. 
$$
Thus, $\Xi_\BB(\ve)$ is an $\ve$-separated set. Moreover, it is easy to see that $\Xi_\BB(\ve)$ is 
maximal. Now, if $x \in \Xi_{\BB}(\ve)$, then $\|x\| \le r_N = \sin \frac12 \t_N = \cos \frac{\pi}{4N}$, 
which implies that $\sqrt{1- \|x\|^2} \ge \sin \frac{\pi}{4N} \ge c \ve$. Thus, \eqref{eq:Xi_B} holds. 
This completes the proof. 
\end{proof}

The unit ball $\BB^d$ equipped with the weight function $(1-\|u\|^2)^{\mu-\f12}$ is also a localizable homogeneous 
space in the sense of \cite{X21}, since it possesses highly localized kernels \cite{PX}. As a result, our construction 
of the maximal $\ve$-separated set on the unit ball also leads to the following proposition, which will not be used
in this paper but is recorded here for possible future use. 
 
 \begin{prop} \label{prop:RemzBall}
Let $W$ be a doubling weight function on $\BB^d$. For $n \in \NN$, let $\chi_{n,\delta}(x,t)$ denote the 
characteristic function of the set $\{x\in \BB^d:  \|x\| \le 1- \frac{\delta}{n^2}\}$. Then for $f \in \Pi_n^d$, 
$1 \le p < \infty$, and every $\delta > 0$, 
\begin{equation}\label{eq:RemzBB}
  \int_{\BB^d} |f(x)|^p W(x) \d x  \le c_\delta \int_{\BB^d} |f(x)|^p \chi_{n,\delta}(x) W(x) \d  x.
\end{equation}
\end{prop}

The proof of the Bernstein inequalities on the cone follows closely the proof on the conic surface. We
will use the maximal function $f_{\b,n}^\ast$ in \eqref{eq:fbn*}, defined similarly by 
\begin{equation} \label{eq:fbn*V}
   f_{\b,n}^\ast (x,t) = \max_{(y,s) \in \VV^{d+1}} \frac{|f(y,s)|}{(1+ n \sd_{\VV}((x,t),(y,s)))^\b}.
\end{equation}
using $ \sd_{\VV}((x,t),(y,s))$ and it remains bounded for any doubling weight $W$ on $\VV^{d+1}$. 
It again satisfies 
\begin{equation} \label{eq:fbn*Vbd}
     \|f\|_{p, W} \le  \|f_{\b,n}^\ast \|_{p, W} \le  c \|f\|_{p, W} 
\end{equation} 
for the doubling weight $W$ on $\VV^{d+1}$ and $1 \le p \le \infty$. 

\begin{lem}
Let $f$ be defined on $\VV^{d+1}$. Define 
\begin{equation} \label{eq:F=f}
      F \big((x,x_{d+1}),t\big) = f(x,t), \qquad \big( (x,x _{d+1}), t\big) \in \VV_0^{d+2}, 
\end{equation}
so that $f$ can be regarded as a function on $\VV_0^{d+2}$. Then
\begin{equation} \label{eq:F*=f*}
    f_{\b,n}^*(x,t) = F_{\b,n}^*\big( (x,x_{d+1}),t\big), \qquad x \in \VV^{d+1}, \quad (x,x_{d+1}) \in \VV_0^{d+2}.
\end{equation}
\end{lem}

\begin{proof}
Since $F$ is independent of $x_{d+1}$, it follows by the definition of $F_{\b,n}^*$ in \eqref{eq:fbn*} 
that $F_{\b,n}^*((x,x_{d+1}), t) = F_{\b,n}^* ((x,-x_{d+1}),t)$. In particular, we can take the maximum
in the definition of $F_{\b,n}^*$ over $\VV_{0,+}^{d+2} = \{(( y,y_{d+1}), t) \in \VV_0^{d+2}: y_{d+1} \ge 0\}$
instead of over $\VV_0^{d+1}$. Then  \eqref{eq:F*=f*} follows from using \eqref{eq:dV=dV0} and replacing 
the maximum in \eqref{eq:fbn*V} as over $(Y,s) \in \VV_{0,+}^{d+2}$.
\end{proof}

\medskip
\begin{proof}[Proof of Theorem \ref{thm:PnCone2}]
Let $f$ be a polynomial of degree at most $n$. Then 
$$
  f(x,t) = \int_{\VV^{d+1} } f(y,s) \Lb_n\left(W_{\g,0}; (x,t), (y,s)\right) W_{\g,0}(y,s) \d y \d s.
$$ 
Let $T_{x,t}$ be a differential operator of $(x,t)$ variables. Then, by \eqref{eq:LnW-Lnw} and 
using the maximal function, we obtain from \eqref{eq:F*=f*}
\begin{align*}
| T_{x,t} f(x,t)| \le f_{\b,n}^*(x,t) & \left[\int_{\VV^{d+1} } \frac{ \left| T_{x,t} \sL_n\left(w_{-1,\g}; (X,t), (Y,s)\right)\right|}
      {\big(1+n \sd_{\VV_0}\big( (X,t),(Y,s)\big) \big)^\b} W_{\g,0}(y,s) \d y \d s \right.\\
        & \left. +   \int_{\VV^{d+1} }  \frac{ \left| T_{x,t} \sL_n\left(w_{-1,\g}; (X,t), (Y_*,s)\right)\right|}
          {\big(1+n \sd_{\VV_0}\big( (X,t),(Y_*,s)\big) \big)^\b} W_{\g,0}(y,s) \d y \d s \right].
\end{align*}
By the identity in Proposition \ref{prop:IntV0V}, we conclude that 
\begin{align}\label{eq:Txtf}
| T_{x,t} f(x,t)| \le f_{\b,n}^*(x,t) \int_{\VV_0^{d+2}}  
   \frac{ \left| T_{x,t} \sL_n\left(w_{-1,\g}; (X,t), (y,s)\right)\right|}
     {\big(1+n \sd_{\VV_0}\big( (X,t),(y,s)\big) \big)^\b} \d \sm(y,s). 
\end{align}
The integral in the right-hand side has been estimated in the proof of Theorems \ref{thm:Bernstein}
and \ref{thm:Bernstein2} in the previous subsection. Hence, using \eqref{eq:fbn*Vbd} and
\eqref{eq:RemzV} when necessary, we have proved the inequalities \eqref{eq:BI-V-1} for $\partial_t$,
\eqref{eq:BI-V-1B} and \eqref{eq:BI-V-1C} for $D_{i,j}$, $1\le i, j \le d$. 

For the $\partial_{x_j}$ we make the following observation. If $G(X) = g(x, \sqrt{t^2-\|x\|^2})$ for 
$x\in \RR^d$ and $|x| \le t$, then
$$
   \sqrt{t^2-\|x\|^2} \partial_j G = \sqrt{t^2-\|x\|^2} \partial_j g - x_j \partial_{d+1} g = D_{j,d+1}^{(X)} G
$$
for $1 \le j \le d$. Consequently, following the procedure in the previous paragraph and using the proof
of Theorem \ref{thm:Bernstein2} for $D_{j,d+1}$, we conclude that
\begin{equation} \label{eq:PhiD-power}
   \left \|(\Phi \partial_{x_j})^\ell f\right\|_{p,W} \le c n^\ell \|f\|_{p,W} \quad \hbox{and} \quad
   \left \| \frac{1}{\sqrt{t}^\ell} (\Phi \partial_{x_j})^\ell f\right\|_{p,W} \le c n^\ell \|f\|_{p,W},
\end{equation}
where the first inequality holds for all $\ell \in \NN$ and the second one holds for $\ell = 1,2$. Furthermore,
the estimate in the Lemma \ref{lem:DkernelV0} shows that 
\begin{align*}
   \frac{\left|\partial_{x_j} \sL_n\left(w_{-1,\g}; (X,t), (y,s)\right)\right|}
   {\big(1+n \sd_{\VV_0}\big( (X,t),(y,s)\big) \big)^\b} 
 &  \le \frac{c_\k n^{d+1} \sqrt{t}  \big(1+n \sd_{\VV_0}\big( (X,t),(y,s)\big) \big)^{-\k+\b}}{ \Phi(x,t)
   \sqrt{\sw_{\g,d}(n;t}\sqrt{\sw_{\g,d}(n;s} } \\
 &  \le \frac{c_\k n^{d+2} \big(1+n \sd_{\VV_0}\big( (X,t),(y,s)\big) \big)^{-\k+\b}}{
   \sqrt{\sw_{\g,d}(n;t}\sqrt{\sw_{\g,d}(n;s} }, 
\end{align*}
if $(x,t)$ satisfies $\Phi(x,t) \ge \delta \sqrt{t} n^{-1}$. Consequently, following the proof of Theorem \ref{thm:Bernstein2},
we conclude that 
$$
  \left | \partial_{x_j} f(x,t) \right | \le  c n^2 f_{\b,n}^*(x,t) \chi_{n,\delta}(x,t), \quad 1 \le j \le d,
$$
where $\chi_{n,\delta}$ is defined in Proposition \ref{prop:RemzV}. Consequently, by \eqref{eq:RemzV} and 
\eqref{eq:fbn*Vbd}, we have proved the inequality $\left \|\partial_{x_j} f \right \|_{p,W} \le c n^2 \|f\|_{p,W}$.
Iterating this inequality proves the first inequality in \eqref{eq:BI-V-2}. Using the identity
$$
  \left( \Phi(x,t) \partial_{x_j}\right)^2 = (\Phi(x,t))^2 \partial_{x_j}^2 - x_j \partial_{x_j},
$$
the first inequality in \eqref{eq:BI-V-2} with $\ell = 1$ and the first inequality in \eqref{eq:PhiD-power}, 
we conclude that 
$$
 \left\| \Phi^2 \partial_{x_j}^2 f \right\|_{p,W} \le  \left\|(\Phi \partial_{x_j})^2 f \right\|_{p,W} + 
  \left\| \partial_{x_j} f \right\|_{p,W} \le c n^2 \|f\|_{p,W},
$$
which proves the second inequality in \eqref{eq:BI-V-2} with $\ell =2$. Since $\Phi^2$ is a polynomial
of $(x,t)$, the second inequality \eqref{eq:BI-V-2} for $\ell > 2$ can be derived by iteration from the case
$\ell =1$ and $2$. Finally, using the second inequality of \eqref{eq:PhiD-power}, the same argument shows that 
$$
 \left\| \frac{1}{t} \Phi^2 \partial_{x_j}^2 f \right\|_{p,W} \le  \left\| \frac{1}{t} (\Phi \partial_{x_j})^2 f \right\|_{p,W} + 
  \left\| \frac{x_j}{t} \partial_{x_j} f \right\|_{p,W} \le c n^2 \|f\|_{p,W},
$$
where we have used $|x_j| \le \|x\| \le t$. This completes the proof. 
\end{proof}

\subsection{Self-adjoint of the spectral operator on the cone}\label{sect:self-adjointV}

We prove Theorem \ref{thm:self-adjointV} that shows the differential operator $\fD_{\mu,\g}$ defined 
in \eqref{eq:V-DE} is self-adjoint. We start with a representation $\fD_{\mu,\g}$ that is of interest in 
itself. It shows, in particular, that the operator $\fD_{\mu,\g}$ can be written in terms of the first-order 
differential operators appeared in Theorem \ref{thm:BI-V}. 

\begin{lem}
The operator $\CD_{\mu,\g}$ can be written as
\begin{align}\label{eq:deq-cone}
  \fD_{\mu,\g} &=  \frac{1}{t^d W_{\mu,\g}(x,t)}    ( \partial_t + t^{-1}\la x,\nabla_x\ra)
      t^{d+1} W_{\mu,\g+1}(x,t) ( \partial_t + t^{-1}\la x,\nabla_x\ra)  \\
        & + t^{-1} \bigg( \sum_{i=1}^d \frac{1}{W_{\mu,\g} (x,t)} \partial_{x_i}  
            \big( W_{\mu+1,\g}(x,t) \partial_{x_i} \big) 
    + \sum_{1 \le i < j \le d } ( D_{i,j}^{(x)})^2 \bigg). \notag
\end{align}
 \end{lem} 

\begin{proof}
We start with an observation that follows from a quick computation,
\begin{align*}
 t (\partial_t + t^{-1} \la x,\nabla_x \ra) (t^2 - \|x\|^2)^{\mu-\f12} & = t (2 \mu-1) (t^2 - \|x\|^2)^{\mu-\f32} (t -t^{-1} \la x,x\ra) \\
    &  =  (2 \mu-1) (t^2 - \|x\|^2)^{\mu-\f12},
\end{align*}
from which it follows that 
\begin{align*}
& \frac{1}{t^d W_{\mu,\g} (x,t)} \left((\partial_t+ t^{-1}\la x,\nabla_x \ra) t^{d+1}
      W_{\mu,\g+1} (x,t) (\partial_t+ t^{-1}\la x,\nabla_x \ra)\right) \\
 & \qquad = t^{-1}(2\mu+d) \la x,\nabla_x \ra + (2\mu+d) \partial_t -(2\mu+d+\g+1)(t \partial_t + \la x ,\nabla \ra)\\
 & \qquad\quad  +t (1-t) (t^{-1} \la x,\nabla_x \ra + \partial_t)^2.
\end{align*}
The last term can be written as  
\begin{align*}
 &t(1-t) (\partial_t+ t^{-1} \la x,\nabla_x \ra )^2  = (1-t) (t \partial_t+ \la x,\nabla_x \ra )(\partial_t+ t^{-1} \la x,\nabla_x \ra ) \\
   & \quad= t (1-t) \partial_t^2+ 2(1-t) \la x,\nabla_x \ra \partial_t - \la x,\nabla \ra^2 + \la x,\nabla_x \ra 
   + t^{-1} \left ( \la x,\nabla \ra^2 - \la x,\nabla_x \ra \right) .
\end{align*}
Putting together shows that 
\begin{align*}
 \fD_{\mu,\g} = & \frac{1}{t^d W_{\mu,\g} (x,t)} \left( (t^{-1} \la x,\nabla_x \ra + \partial_t) t^{d+1}
       W_{\mu,\g+1}(x,t) (t^{-1}\la x,\nabla_x \ra + \partial_t)\right) \\
    & +t^{-1} \left (t^2 \Delta_x - \la x,\nabla_x\ra^2 - (2\mu + d-1) \la x,\nabla_x\ra\right).
\end{align*}
To conclude the proof, we need the following identity, 
\begin{align*}
 & t^2 \Delta_x - \la x,\nabla_x\ra^2 - (2\mu + d-1) \la x,\nabla_x\ra \\
 & \quad =    \frac{1}{W_{\mu,\g} (x,t)} \sum_{i=1}^d  \partial_{x_i}  \big( W_{\mu+1,\g}(x,t) \partial_{x_i} \big) 
    + \sum_{1 \le i < j \le d } ( D_{i,j}^{(x)})^2,
\end{align*}
which we claim to hold. Indeed, if we dilate this identity by $x = t y$, and use the relations $ t^2 \Delta_x=\Delta_y$,
$ \la x,\nabla_x \ra= \la y,\nabla_y\ra$, $  D_{i,j}^{(x)} = D_{i,j}^{(y)}$, and 
\begin{equation} \label{eq:x=ytBall}
  \frac{1}{W_{\mu,\g}(x,t)}    \partial_{x_i}  \left(  W_{\mu+1,\g}(x,t) \partial_{x_i} \right)
  =  \frac{1}{\varpi_\mu (y)}  \partial_{y_i}  \left( \varpi_{\mu+1}(y) \partial_{y_i} \right), 
\end{equation}
where $\varpi_\mu(y) = (1-\|y\|^2)^{\mu-\f12}$, we see that the claimed identity becomes 
\begin{align*}
 & \Delta_y - \la y,\nabla_y\ra^2 - (2\mu + d-1) \la y,\nabla_y\ra \\
   & \qquad =  \frac{1}{\varpi_\mu (y)} \sum_{i=1}^d  \partial_{y_i}  \left(  \varpi_{\mu+1}(y) \partial_{y_i} \right) 
    + \sum_{1 \le i < j \le d } ( D_{i,j}^{(y)})^2,
\end{align*}
which can be easily verified and is in fact the decomposition of the spectral differential operator for the 
unit ball $\BB^d$ (cf. \cite[Section 5.2]{DX}). This completes the proof.
\end{proof}

\noindent
\begin{proof}[Proof of Theorem  \ref{thm:self-adjointV}]
Using the decomposition \eqref{eq:deq-cone} of $\fD_{\mu,\g}$ and making a change of variable
$x = yt$, we can write
\begin{align*}
 &  \int_{\VV^{d+1}} \fD_{\mu,\g} \cdot g W_{\mu,\g} \d x\d t 
  =  \int_0^1 t^d \int_{\BB^d} \fD_{\mu,\g} f(t y,t) \cdot g(t y,t) W_{\mu,\g}(t y,t) \d y\d t \\
& =   \int_{\BB^d}  \int_0^1 ( \partial_t + t^{-1}\la x,\nabla_x\ra)
      \left( t^{d+1} W_{\mu,\g+1}(x,t) ( \partial_t + t^{-1}\la x,\nabla_x\ra) \right) f(t y,t) \cdot g(x,t)\d t  \d y \\
& \qquad\quad + \int_0^1 t^{2\mu+d-2}   \sum_{i=1}^d  \int_{\BB^d} \partial_{y_i}  
    \left( (1-\|y\|^2)^{\mu+\f12} \partial_{y_i} \right) 
   f(t y,t) \cdot g(x,t)   \d y (1-t)^\g \d t   \\
& \qquad\quad  + \sum_{1 \le i < j \le d }  \int_0^1 t^{2\mu+ d-2} \int_{\BB^d} ( D_{i,j}^{(y)})^2 
     f(t y,t) \cdot g(t y,t) \varpi_\mu (y) \d y\, (1-t)^\g\d t,    
\end{align*}
where we have used \eqref{eq:x=ytBall} in the second term in the righthand side and 
$  D_{i,j}^{(x)} =  D_{i,j}^{(y)}$ in the third term in the righthand side. Assuming $\g > 0$, we observe that 
$$
   \frac{\d}{\d t} [f(t y,t)] = \big( (t^{-1} \la x,\nabla \ra + \partial_t) f \big) (t y, t), 
$$ 
so that we can do integration by parts with respect to the variable $t$ in the first integral in the right-hand side.  
Integration by parts also shows that the second integral in the right-hand side is equal to 
\begin{align*}
  & - \int_0^1 t^{2\mu+d-2}   \sum_{i=1}^d \int_{\BB^d}  \partial_{y_i} f(t y,t) \cdot  
     \partial_{y_i} g(t y,t)\varpi_{\mu+1}(y) \d y (1-t)^\g\d t \\
  &\qquad  =   - \int_0^1 t^{d-1} \sum_{i=1}^d \int_{\BB^d}     \partial_{x_i} f(x, t) \cdot  
     \partial_{x_i} g(x,t)W_{\mu+1,\g}(x,t) \d x  \d t, 
\end{align*}
which is the second integral in the right-hand side of \eqref{eq:intJacobiSolid}. Next, sing the result on
the unit ball, the third integral in the right-hand side is equal to 
\begin{align*}
 & - \sum_{1 \le i < j \le d }  \int_0^1 t^{2\mu+ d-2} \int_{\BB^d} D_{i,j}^{(y)} 
    f(t y,t) \cdot D_{i,j}^{(y)} g(t y,t) \varpi_\mu (y) \d y\, (1-t)^\g\d t \\
 & \qquad \quad = - \sum_{1 \le i < j \le d }  \int_0^1 t^{d-1} \int_{\BB^d} D_{i,j}^{(x)} 
    f(x,t) \cdot D_{i,j}^{(x)} g(x,t) W_{\mu,\g} (x,t) \d x\,\d t,
    \end{align*}
which gives the third integral in the right-hand side of \eqref{eq:intJacobiSolid}. Finally, analytic continuation 
shows that the identity holds for $\g > -1$. 
\end{proof}

\end{document}